%\magnification=\magstep1

\input amssym.def
\input amssym.tex
\input xy
\xyoption{all}
\CompileMatrices
\tracingstats=2

%%%%%%       M A C R O S       %%%%%%%

%%%%%%    FONTS    %%%%%%%

\font\twbf=cmbx12
\font\cyr=wncyr10

\font\tenfun=rsfs10
\font\sevenfun=rsfs7
\font\fivefun=rsfs5

\newfam\funfam
\def\fun{\fam\funfam\tenfun}
\textfont\funfam=\tenfun
\scriptfont\funfam=\sevenfun
\scriptscriptfont\funfam=\fivefun

%%%%%%%%%%%%%%%%%%%%%%%%%%%%%%%%%%%%%%%%%%

\def\bn{\bigskip\noindent}
\def\sn{\smallskip\noindent}

\def\({\left(}
\def\){\right)}

\def\ddt{{\scriptstyle\scriptstyle\bullet}}

\def\o{\overline}

\def\C{{\bf C}}
\def\F{{\bf F}}
\def\K{{\cal K}}
\def\N{{\bf N}}
\def\O{{\cal O}}
\def\wp{{\goth p}}
\def\Q{{\bf Q}}
\def\R{{\bf R}}

\def\Z{{\bf Z}}
\def\sha{\hbox{\cyr X}}

\def\qb{\overline{\bf Q}}
\def\qp{{\bf Q}_p}

\def\zp{{\bf Z}_p}

\def\Bd{B^t}
\def\ki{K_\infty}

\def\coh #1#2#3{H^{#1}(#2,#3)}
\def\cohur #1#2{H^1_{ur}(#1,#2)}
\def\cohf #1#2#3{\widetilde H^{#1}_f(#2,#3)}
\def\cohbk #1#2{H^1_f(#1,#2)}
\def\cohiwf #1#2#3#4{\widetilde H^{#1}_{f,{\rm Iw}}(#2/#3,#4)}
\def\Cff #1#2{\widetilde C^\ddt_f({#1,#2})}

\def\Tam #1#2#3{{\rm Tam}_{#1}(#2,#3)}

\def\fr #1#2{{{#1} \over {#2}}}

\def\ad{{\rm ad}}
\def\Coker{{\rm Coker}}
\def\Cone{{\rm Cone}}
\def\cont{{\rm cont}}
\def\cor{{\rm cor}}
\def\cork{{\rm cork}}
\def\divv{{\rm div}}
\def\End{{\rm End}}
\def\Ext{{\rm Ext}}

\def\gen{{\rm gen}}
\def\Gal{{\rm Gal}}
\def\Hom{{\rm Hom}}
\def\Imm{{\rm Im}}
\def\id{{\rm id}}
\def\Ker{{\rm Ker}}
\def\ord{{\rm ord}}
\def\Pic{{\rm Pic}}
\def\res{{\rm res}}
\def\rk{{\rm rk}}
\def\Sel{{\rm Sel}}
\def\sgn{{\rm sgn}}
\def\str{{\rm str}}
\def\tors{{\rm tors}}
\def\Tr{{\rm Tr}}

\def\ho{\hookrightarrow}
\def\iso{\buildrel \sim \over \longrightarrow}
\def\lo{\longrightarrow}
\def\Lo{\Longrightarrow}

\def\dirlim{\mathop{\vtop{\ialign{##\crcr$\hfill{\lim}\hfil$\crcr
\noalign{\kern1pt\nointerlineskip}\rightarrowfill\crcr\noalign
{\kern -2pt}}}}\limits}

\def\invlim{\mathop{\vtop{\ialign{##\crcr$\hfill{\lim}\hfil$\crcr
\noalign{\kern1pt\nointerlineskip}\leftarrowfill\crcr\noalign
{\kern -2pt}}}}\limits}

\def\normalbaselines{\baselineskip20pt\lineskip3pt\lineskiplimit3pt}
\def\mapr#1{\smash{\mathop{\longrightarrow}\limits^{#1}}}
\def\lomapr#1{\smash{\mathop{\relbar\joinrel\longrightarrow}\limits^{#1}}}
\def\mapd#1{\Big\downarrow\rlap{$\vcenter{\hbox{$\scriptstyle#1$}}$}}

\def\hb{\hfill\break}

\centerline{\twbf Kolyvagin's result on the vanishing of $\sha(E/K)[p^\infty]$
and its consequences}
\vskip4pt
\centerline {\twbf for anticyclotomic Iwasawa theory}

\vskip10pt
\centerline{\bf Ahmed Matar, Jan Nekov\'a\v{r}}

\bn
Abstract: We discuss improvements of Kolyvagin's classical result about
the vanishing of the $p$-primary part of the Tate--\v{S}afarevi\v{c} group
of an elliptic curve $E$ (defined over ${\bf Q}$) over an imaginary quadratic
field $K$ satisfying the Heegner hypothesis for which the basic Heegner point
$y_K \in E(K)$ is not divisible by an odd prime $p$. Combining Kolyvagin's
theorem with a new abstract Iwasawa-theoretical result, we deduce, under
suitable assumptions, that similar vanishing holds for all layers in the
anticyclotomic ${\bf Z}_p$-extension of $K$.

\vskip16pt

\centerline{\bf 0. Introduction}

\bn
{\bf (0.1)} \ Let $E$ be an elliptic curve over $\Q$ of conductor $N$ and $K$
an imaginary quadratic field of discriminant $D_K$ in which all primes
dividing $N$ split. Fix a modular parameterisation $\varphi : X_0(N) \lo E$ and
an ideal ${\cal N} \subset O_K$ such that $\O_K/{\cal N} \simeq \Z/N\Z$.
The basic Heegner point $y_K \in E(K)$ attached to these data is, by
defi\-nition, the trace $y_K := \Tr_{H_1/K}(y_1)$ of the Heegner point
of conductor one $y_1 := \varphi([\C/O_K \lo \C/{\cal N}^{-1}]) \in E(H_1)$
defined over the Hilbert class field $H_1$ of $K$.

\sn
{\bf (0.2)} \ If $y_K \not\in E(K)_{\tors}$ and $D_K \not= -3, -4$, Kolyvagin
[K1, Thm. A] proved that the groups $E(K)/\Z y_K$ and $\sha(E/K)$ are finite,
and that $\# \sha(E/K)$ divides $[E(K) : \Z y_K]^2$ multiplied by a product
of several error terms. The $p$-primary parts of these error terms vanish
in the following situation (each of the respective assumptions (a), (b) and (c)
implies that the corresponding error term $a, b, c$ in [K1, Cor. 11, Cor. 12,
Cor. 13] is relatively prime to $p$; the error term $d$ is equal to $1$, since
$p \not= 2$).

\proclaim{(0.3) Theorem (Kolyvagin, special case of [K1, Cor. 13])}. Assume
that $D_K \not= -3, -4$ and that $p \not= 2$ is a prime number satisfying
the following conditions.\hb
(a) $\forall n_1, n_2\geq 0 \quad H^1(K(E[p^{n_1+n_2}])/K, E[p^{n_1}]) = 0$.\hb
(b) Neither of the $(\pm 1)$-eigenspaces $E[p]^{\pm}$ for the action of complex
conjugation is stable under the action of $G_{\Q} :=\Gal(\qb/\Q)$. Equivalently,
the ${\rm (mod} \; p{\rm )}$ Galois representation $\o\rho_{E,p} : G_{\Q} \lo
{\rm Aut}_{\F_p}(E[p]) \simeq GL_2(\F_p)$ is irreducible.\hb
(c) $E(K)[p] = 0$.\hb
If $y_K \not\in E(K)_{\tors}$, then $E(K)/\Z y_K$ is finite and
$$
p^{m_0} \sha(E/K)[p^\infty] = 0,\qquad \# \sha(E/K)[p^\infty]
{\rm\ divides\ } p^{2 m_0},
$$
where $m_0 := \sup\{m\geq 0 \mid y_K \in p^m E(K)\}$ (thus $E(K) \otimes\zp
\simeq \zp$ and $p^{m_0} = [E(K) \otimes\zp : \zp(y_K \otimes 1)]$).

\sn
{\bf (0.4)} \ For $p \not= 2$, the assumption (b) in Theorem 0.3 implies (c).
Moreover, the assumptions (a), (b) and (c) are satisfied
%,in particular,
if $\o\rho_{E,p}$ has ``big image" (e.g., if it is surjective).

Gross [G] gave a self-contained account of Kolyvagin's proof of Theorem 0.3
in the simplest case when $\o\rho_{E,p}$ is surjective and $m_0 = 0$. One step
in the argument ([G, beginning of \S{9}]) required an additional assumption
$p \nmid D_K$.

\proclaim{(0.5) Theorem ([G, Prop. 2.1, Prop. 2.3])}. Assume that $D_K \not=
-3, -4$ and that $p \nmid 2 D_K$ is a prime number  for which $\o\rho_{E,p} :
G_{\Q} \lo GL_2(\F_p)$ is surjective. If $y_K \not\in p E(K)$, then $E(K)
\otimes \zp = \zp y_K \simeq \zp$ and $\sha(E/K)[p^\infty] = 0$.

\sn
{\bf (0.6)} \ In [K2], Kolyvagin proved the following structure theorem for
$\sha(E/K)[p^\infty]$, which refines Theorem 0.3 (under the ``big image"
assumption for the $p$-adic Galois representation $\rho_{E,p} : G_{\Q} \lo 
{\rm Aut}_{\zp}(T_p(E)) \simeq GL_2(\zp)$).

\proclaim{(0.7) Theorem (Kolyvagin, [K2, Thm. C, Thm. D])}. Assume that
 $D_K \not= -3, -4$ and that $p \not= 2$ is a prime number for which
$\rho_{E,p} : G_{\Q} \lo GL_2(\zp)$ has ``big image" (e.g., that $\rho_{E,p}$
is surjective). If $y_K \not\in E(K)_{\tors}$, then
$$
\sha(E/K)[p^\infty] \simeq X \oplus X,\qquad X\simeq \bigoplus_{i\geq 0}
\Z/p^{m_i - m_{i+1}}\Z,\qquad m_0 \geq m_1 \geq \cdots \geq m_\infty :=
\inf m_i,
$$
where $m_0$ is as in Theorem 0.3 and $m_i$ for $i > 0$ is defined in a similar
way in terms of certain linear combinations of Heegner points of higher
conductors. In particular, $\# \sha(E/K)[p^\infty] = p^{2 (m_0 - m_\infty)}$.

\sn
{\bf (0.8)} \ The divisibility $\# X \mid p^{m_0}$ was reproved by Howard
[H1, Thm. A] using the formalism of anti\-cyclotomic Kolyvagin systems,
under the assumptions that $\rho_{E,p}$ is surjective, $D_K \not= -3, -4$
and $p\nmid  2 N D_K$.

\sn
{\bf (0.9)} \ For $p\not= 2$, the condition (a) in Theorem 0.3 was studied in
detail by Cha [Ch, Thm. 2], who showed that it is satisfied if $p\nmid D_K$,
$p^2 \nmid N$ and $E(K)[p] = 0$, except when $p = 3$ and $\o\rho_{E,3}(G_K) =
\pmatrix{\F_3^\times & \F_3\cr 0 & 1\cr}$.
Therefore the conclusions of Theorems 0.3 and 0.5 hold (for $D_K \not= -3, -4$)
whenever $p\nmid 2 D_K$, $p^2 \nmid N$ and $\o\rho_{E,p}$ is irreducible. He
also showed [Ch, Thm. 21, Rmk. 25] that the statement of Theorem 0.7 holds
under the same assumptions.

\sn
{\bf (0.10)} \ The authors of a collective article [GJPST] had made an attempt
at generalising Cha's results. However, the cohomological calculations in
[GJPST, Lemma 5.7, Lemma 5.9] and [GJPST, proof of Proposition 5.4] are
incorrect (see [LW, Lemma 8]), the statement of [GJPST, Proposition 5.8]
is correct but the proof is not, and the discussion of Kolyvagin's
method (in the form presented in [G]) in [GJPST, \S{5}] is seriously flawed.
In particular, the assertion to the effect that the surjectivity of
$\o\rho_{E,p}$ in Theorem 0.5 can be replaced by the vanishing of the
groups $H^i(K(E[p])/K, E[p])$ for $i = 1, 2$ and of $E^\prime(K)[p]$
for all $\Q$-isogenies $E \lo E^\prime$, is incorrect, for the following
reason: the current state of the art requires an irreducibility assumption
for $\o\rho_{E,p}$ (or its restriction to $G_K$) in order to obtain, by
Kolyvagin's method, an upper bound on the size of $\sha(E/K)[p^\infty]$
without any error terms. As a result, [GJPST, Thm. 3.7] remains unproved.

\sn
{\bf (0.11)} \ Lawson and Wuthrich [LW] extended and simplified the
cohomological calculations of [Ch], and corrected various mistakes
from [GJPST]. In [LW, Thm 1, Thm 2], they gave a complete classification of
pairs $(E, p)$ consisting of an elliptic curve $E$ over $\Q$ and a prime
number $p \not= 2$ for which $H^1(\Q(E[p])/\Q, E[p]) \not= 0$ (and similarly
for $H^1(\Q(E[p^n])/\Q, E[p^n]) \not= 0$, where $n > 1$ and $p > 3$).
They also classified pairs $(E, p)$ for which $H^2(\Q(E[p])/\Q, E[p]) \not= 0$.

Their results imply that the condition (a) in Theorem 0.3 (for $p\not= 2$)
is always satisfied if $\o\rho_{E,p}$ is irreducible. Consequently, the
conclusions of Theorems 0.3 and 0.7 hold (for $D_K \not= -3, -4$) if
$\o\rho_{E,p}$ is irreducible and $p\not= 2$.

However, the claims made in [LW, Thm. 14] about the validity of Theorem 0.3
in situations when (a) holds but $\o\rho_{E,p}$ is reducible are unjustified,
for reasons explained in 0.10.

We recall the methods of [Ch] and [LW] and prove a mild generalisation
of some of their results in \S{5}. We also prove the following variant
of Theorem 0.5.

\proclaim{(0.12) Theorem (= Theorem 6.7)}. Assume that $p \not= 2$ and that
$E[p]$ is an irreducible $\F_p[G_{\Q}]$-module (which implies that $E(K)[p]
= 0$).\hb
(1) Assume that $(K, p) \not= (\Q(\sqrt{-3}), 3)$. If $y_K \not\in p E(K)$,
then
$$
E(K) \otimes\zp = \zp(y_K \otimes 1) \simeq \zp,\qquad \sha(E/K)[p^\infty] = 0.
$$
(2) Assume that $(K, p) = (\Q(\sqrt{-3}), 3)$; then $y_K \in 3 E(K)$.
If $3\nmid a_3$, assume, in addition, that $\rho(G_K)$ is not a cyclic group
of order four. If $y_K \not\in 3^2 E(K)$, then
$$
\Z_3 \simeq E(K) \otimes\Z_3 \supset 3 E(K) \otimes\Z_3 = \Z_3 (y_K \otimes 1),
\qquad \sha(E/K)[3^\infty] = 0.
$$

\sn
{\bf (0.13)} \ We now turn to Iwasawa-theoretical results. Fix a prime number
$p$ and denote by $\ki = \bigcup_{n\geq 1} K_n$ the anticyclotomic
$\zp$-extension of $K$. In this case $\Gamma := \Gal(\ki/K) \simeq \zp$, $K_n =
K_\infty^{\Gamma_n}$, where $\Gamma_n = \Gamma^{p^n} \simeq p^n \zp$, and
$\Gal(\ki/\Q) = \Gamma \rtimes \{1, c\}$, with complex conjugation $c$ acting
on $\Gamma$ by $g \mapsto g^{-1}$. Denote by $\Lambda := \zp[[\Gamma]]$ the
Iwasawa algebra of $\Gamma$.

\sn
{\bf (0.14)} \ From now on until the end of Introduction we assume that
$p \not= 2$ and that $E$ has good ordinary reduction at $p$. The Selmer
module $\Sel_{p^\infty}(E/\ki) := \dirlim_n \Sel_{p^\infty}(E/K_n)$ (resp.
$S_p(E/\ki) := \invlim_n S_p(E/K_n)$) (see \S{1.4} for the notation) is
a $\Lambda$-module of cofinite (resp. finite) type, of corank (resp. rank)
equal to one, as predicted by one of Mazur's conjectures formulated in
[Mz, \S{18}]. This conjecture is a consequence of another conjecture of Mazur
[Mz, \S{19}] (proved independently by Cornut [Co1, Co2] and Vatsal [Va])
combined with an Euler system argument along the tower $\ki/K$ ([B, Thm. A]
under some additional assumptions; the general case is proved in [N1, \S{2}]
together with [N3, Thm. 3.2]; see also [H1, Thm. B]). Another proof
of [B, Thm. A], which had applications to the study of the anticyclotomic
$\mu$-invariant, was given in [M, Thm. A].

In [B, Thm. B], Bertolini also proved a $\Lambda$-adic variant of Kolyvagin's 
annihilation result [K1, Cor. 12] for the torsion submodule of the Pontryagin
dual of $\Sel_{p^\infty}(E/\ki)$ (assuming the validity of Mazur's conjecture
[Mz, \S{19}]). This result was subsequently generalised by Howard [H1, Thm. B],
who proved one half of a conjecture of Perrin-Riou [PR, Conj. B] for Heegner
points along $\ki/K$, namely, a $\Lambda$-adic variant of Kolyvagin's result
$\# X \mid p^{m_0}$ (in the notation of Theorem 0.7).

\sn
{\bf (0.15)} \ The proofs of [B, Thm. B] and [H1, Thm. B] relied on fairly
detailed arguments involving the Euler system and the Kolyvagin system
of Heegner points along $\ki/K$, respectively. The main insight of the present
work is that in the simplest case when $y_K \not\in p E(K)$, one can obtain
--- under certain assumptions --- precise information about the structure of
the $\zp[\Gamma/\Gamma_n]$-modules $E(K_n) \otimes \qp/\zp$ and
$\sha(E/K_n)[p^\infty]$ from Theorem 0.5 (and its variant Theorem 0.12) by
purely Iwasawa-theoretical
methods, combined with the norm relations for the Heegner points of $p$-power
conductor, without applying any Euler system arguments along the tower
$\ki/K$. The following results are proved in \S{4}.

\proclaim{(0.16) Theorem (= Theorem 4.8)}. If $p \not= 2$ is a prime number
such that\hb
(a) $E(K)[p] = 0$,\hb
(b) $p\nmid N \cdot a_p \cdot (a_p - 1) \cdot c_{\rm Tam}(E/\Q)$,\hb
(c) $y_K \not\in E(K)_{\tors}$,\hb
(d) $\rk_{\Z}\, E(K) = 1$ and $\sha(E/K)[p^\infty] = 0$,\hb
then $\sha(E/\ki)[p^\infty] = 0$ and the Pontryagin dual of $E(\ki) \otimes
\qp/\zp = \Sel_{p^\infty}(E/\ki)$ is a free module of rank one over
$\zp[[\Gal(\ki/K)]]$.

\proclaim{(0.17) Theorem (= Theorem 4.9)}. If $p \not= 2$ is a prime number
such that\hb
(a) $E(K)[p] = 0$,\hb
(b') $p\nmid N \cdot a_p \cdot (a_p - 1) \cdot (a_p - \eta_K(p))\cdot
c_{\rm Tam}(E/\Q)$,\hb
(c') $y_K \not\in p E(K)$,\hb
(d) $\rk_{\Z}\, E(K) = 1$ and $\sha(E/K)[p^\infty] = 0$,\hb
then, for every intermediate field $K \subset L \subset \ki$,
$\sha(E/L)[p^\infty] = 0$ and the Pontryagin dual of $E(L) \otimes \qp/\zp =
\Sel_{p^\infty}(E/L)$ is a free module of rank one over
$\zp[[\Gal(L/K)]]$. For every integer $n\geq 0$, $\rk_{\Z}\, E(K_n) = p^n$,
$\sha(E/K_n)[p^\infty] = 0$ and $E(K_n) \otimes \zp$ is generated over
$\zp[\Gal(K_n/K)]$ by the traces to $K_n$ of the Heegner points of $p$-power
conductor.

\sn
{\bf (0.18)} \ Above, $a_p$ denotes the $p$-th coefficient of the $L$-function
$L(E/\Q, s) = \sum_{n\geq 1} a_n n^{-s}$, the value $\eta_K(p)$ is equal to
$1, -1, 0$, respectively, if $p$ splits, is inert, or is ramified in $K/\Q$,
and $c_{\rm Tam}(E/\Q) = \prod_{\ell \mid N} c_{\rm Tam, \ell}(E/\Q)$ is the
product of the local Tamagawa factors of $E$ at all primes of bad reduction.

\sn
{\bf (0.19)} \ If $K = \Q(\sqrt{-3})$ and $p = 3$, the conditions (a) and
(c') in Theorem 0.17 cannot be satisfied simultaneously, by Proposition~4.11
below.  In general, (a) and (c') should imply both (d) and
$p\nmid c_{\rm Tam}(E/\Q)$ (see (6.2.1)).

\sn
{\bf (0.20)} \ What is the role of the individual assumptions in Theorem 0.16
and Theorem 0.17? The condition (a) implies that $E(\ki)[p]= 0$. The assumption
$p \nmid N \cdot a_p$ is equivalent to $E$ having good ordinary reduction
at $p$, and the remaining part $p \nmid (a_p - 1) \cdot c_{\rm Tam}(E/\Q)$
of (b) ensures (when combined with (a)) that Mazur's control theorem holds
along the tower $\ki/K$ without any error terms: $\Sel_{p^k}(E/K_n) \iso
\Sel_{p^k}(E/\ki)^{\Gamma_n}$ for all $k, n\geq 0$. The condition (d) implies
that $\Sel_{p^\infty}(E/K) \simeq \qp/\zp$. Finally, the norm relations
for the Heegner points of $p$-power conductor imply that, for a suitable
non-zero element $m\in \zp$, the multiple $y_K \otimes m \in E(K) \otimes \zp$
is a universal norm from the projective system $\{E(K_n) \otimes \zp\}$, and
the condition $p\nmid (a_p - 1) \cdot (a_p - \eta_K(p))$ ensures that
$m \in {\bf Z}_p^\times$.

\sn
{\bf (0.21)} \ One can combine Theorem 0.17 with the Euler system results
over $K$ (but not over $\ki$) discussed in 0.1--0.11. Kolyvagin's result
alluded to in 0.2 tells us that the condition $\rk_{\Z}\, E(K) = 1$ in
Theorem 0.16(d) follows from (c), and therefore can be dropped. Likewise,
the condition (d) in Theorem 0.17 follows from (c'), whenever the conclusions
of Theorem 0.5 hold. Combining Theorem 0.5 (with weaker assumptions, supplied
by [Ch], [LW] and Theorem 6.7(1)) with Theorem 0.17, we obtain the following
result.

\proclaim{(0.22) Theorem (= Theorem 6.9)}. If $p \not= 2$ is a prime number
such that $E[p]$ is an irreducible $\F_p[G_{\Q}]$-module, $p\nmid N \cdot
a_p \cdot (a_p - 1) \cdot (a_p - \eta_K(p))\cdot c_{\rm Tam}(E/\Q)$ and
$y_K \not\in p E(K)$, then the conclusions of Theorem 0.17 hold.

\sn
{\bf (0.23)} \ The case $K = \Q(\sqrt{-3})$, $p = 3$ is different,
as already mentioned in Theorem 0.12 and in 0.19. The point is that,
if $E(K)[3] = 0$, then $y_K = 3 z_K$, where $z_K \in E(K)$ is a linear
combination of the traces to $K$ of the Heegner points of conductors $1$
and $q$, for any prime $q \nmid 3N$ satisfying $a_q \not\equiv 1 + 
\eta_K(q) \; ({\rm mod}\ 3)$ (there are infinitely many such primes $q$).

\proclaim{(0.24) Theorem (= Theorem 6.10)}.
Assume that $K = \Q(\sqrt{-3})$, $p = 3$, $E[3]$ is an irreducible
$\F_3[G_{\Q}]$-module, $\rho(G_K)$ is not a cyclic group of order four
and $3\nmid a_3 \cdot (a_3 - 1) \cdot c_{\rm Tam}(E/\Q)$.
If $y_K \not\in 3^2 E(K)$, then the conclusions of Theorem 0.17 hold,
with the following modification: each $E(K_n) \otimes {\bf Z}_3$ is
generated over ${\bf Z}_3[\Gal(K_n/K)]$ by the traces to $\ki$ of the
Heegner points of conductors dividing $3^\infty q$, for any prime $q$
as in 0.23.

\sn
{\bf (0.25)} \ Analogous results hold for anticyclotomic
${\bf Z}_p^m$-extensions and basic CM points on abelian varieties of
$GL(2)$-type with real multiplication occurring as simple quotients of
Jacobians of Shimura curves over totally real number fields. This will
be discussed in a separate publication.

\sn
{\bf (0.26)} \ Let us describe the contents of this article in more detail.
The goal of \S{1}-\S{3} is to prove two abstract results (Theorems~3.4
and 3.5) on Selmer groups of $\wp$-ordinary abelian varieties in dihedral
Iwasawa theory. The framework is general enough to apply in the context
of 0.25, not just in the situation involving classical Heegner points
on elliptic curves. In \S{4} we recall the norm relations for Heegner
points and combine them with Theorems~3.4 and 3.5 in order to deduce
Theorems~0.16 and 0.17. In \S{5}-\S{6} we give a proof of Kolyvagin's
result on vanishing of $\sha(E/K)[p^\infty]$ in the form of Theorem~0.12.
When combined with Theorems~0.16 and 0.17, this implies Theorems~0.22 and
0.24. Again, the general theory developed in \S{5} is applicable in the
context of 0.25.

\sn
{\bf Acknowledgements.} \ Some of the work on this article was carried out by
the second named author when he was visiting Centre Interfacultaire Bernoulli
(CIB) at Ecole Polytechnique F\'ed\'erale de Lausanne during the semester
``Euler systems and special values of $L$-functions" in fall 2017, and when he
was staying at Imperial College London as an ICL--CNRS fellow in spring 2018.
He is grateful to both institutions for their generous support.

\vskip16pt

\centerline{\bf 1. Generalities}

\bn
{\bf (1.1)} \ Throughout \S{1}-\S{3},

\medskip
\item{$\bullet$} \ for any perfect field $k$, denote by $G_k = \Gal(\o k/k)$
its absolute Galois group.
\item{$\bullet$} \ For an integer $n\geq 1$ invertible in $k$, denote by
$\chi_{n,k} : G_k \lo (\Z/n\Z)^\times$ the cyclotomic character given
by the action of $G_k$ on $\mu_n(\o k)$.
\item{$\bullet$} \ $K$ is a number field.
\item{$\bullet$} \ $Fr(v)$ will always denote the arithmetic Frobenius
element.
\item{$\bullet$} \ $p$ is a prime number; if $K$ is not totally imaginary,
we assume that $p \not= 2$.
\item{$\bullet$} \ $B$ is an abelian variety over $K$ with {\bf good reduction}
at all primes of $K$ above $p$; let $\Bd$ be the dual abelian variety.
\item{$\bullet$} \ If $v$ is a finite prime of a finite extension $L$ of $K$,
denote by $B_v$ the N\'eron model of $B \otimes_K L_v$ over $O_{L_v}$, by
$\widetilde B_v$ its special fibre (over the residue field $k(v)$ of $v$),
and by $\pi_0(\widetilde B_v) = \widetilde B_v/ \widetilde B_v^\circ$ the
$G_{k(v)}$-module of its connected components.
\item{$\bullet$} \ $M$ is a totally real number field with ring of integers
$O_M$.
\item{$\bullet$} \ We are given a ring morphism $i : O_M \lo \End(B)$
and an $O_M$-linear isogeny $\lambda : B \lo \Bd$ which is symmetric
in the sense that $\lambda = \lambda^t$. Above, we use a scheme-theoretic
notation: the ring of endomorphisms of $B$ defined over a field $L$
containing $K$ is denoted by $\End(B \otimes_K L)$ (not by $\End_L(E)$).

\medskip
Throughout, one can replace $O_M$ by any order in $M$ whose index in $O_M$
is prime to $p$, but the current setting is sufficient for the arithmetic
applications we have in mind.

\sn
{\bf (1.2)} \ The decomposition

$$
O_M \otimes \zp = \prod_{\wp\mid p} O_{M_\wp}
$$
(where $\wp$ runs through all primes of $M$ above $p$) induces decompositions

$$
B[p^\infty] = \bigoplus_{\wp} B[\wp^\infty],\qquad
T_p(B) = \bigoplus_{\wp} T_{\wp}(B).
$$
Fix, once for all, a prime $\wp \mid p$ in $M$ and set

$$
\O := O_{M_\wp},\qquad \K := M_\wp,\qquad A:= B[\wp^\infty],\qquad
T:= T_{\wp}(B).
$$
Throughout \S{1}-\S{3}, we assume that

\medskip
\item{$\bullet$} \ $B$ has good {\bf $\wp$-ordinary reduction} at each
prime $v$ of $K$ above $p$ in the sense that
$$
\rk_{\O} \, T_{\wp}(\widetilde B_v) = {1\over 2}\, \rk_{\O} \, T_{\wp}(B)
\qquad (= \dim(B)/[M:\Q]).
$$
This condition is weaker than requiring $B$ to have good ordinary reduction
at $v$ (which is equivalent to $B$ having good ${\goth p}^\prime$-ordinary
reduction at $v$ for all ${\goth p}^\prime \mid p$ in $M$).

\sn
{\bf (1.3) Pontryagin duality.} \ For any discrete or compact topological
$\zp$-module $X$, let us denote by

$$
D(X) := \Hom_{\cont,\zp}(X, \qp/\zp)
$$
the Pontryagin dual of $X$. In the special case when $X$ is a topological
$\O$-module, so is $D(X)$, and there are canonical isomorphisms of $\O$-modules

$$
\displaylines{
D(X) \iso \Hom_{\cont,\O}(X, \Hom_{\zp}(\O, \qp/\zp)),\cr
\Hom_{\zp}(\O, \qp/\zp) \iso \Hom_{\zp}(\O, \zp) \otimes_{\zp}
\qp/\zp = \Hom_{\zp}(\O, \zp) \otimes_{\O} \K/\O,\cr}
$$
where $\Hom_{\zp}(\O, \zp)$ is a free $\O$-module of rank one.
A choice of an isomorphism of $\O$-modules

$$
\O \iso\Hom_{\zp}(\O, \zp)
\eqno (1.3.1)
$$
is equivalent to choosing a generator $a \in {\fun D}_{\O/\zp}^{-1}$
of the inverse different, via the pairing

$$
\O \times \O \lo \zp,\qquad (x, y) \mapsto \Tr_{\K/\qp}(axy).
\eqno (1.3.2)
$$
As in [N2, (0.4.1)], we let

$$
T^* := D(A),\qquad A^* := D(T).
$$
The Weil pairing

$$
(\; ,\;) : T_p(B) \times T_p(\Bd) \lo \zp(1)
$$
is $\zp$-bilinear and $G_K$-equivariant. It satisfies $(\alpha x, y) =
(x, \alpha^t y)$, for all $\alpha \in \End(B)$ (where $\alpha^t$ denotes
the dual isogeny to $\alpha$). In particular, it induces an eponymous pairing

$$
(\; ,\;) : T_{\wp}(B) \times T_{\wp}(\Bd) \lo \zp(1)
\eqno (1.3.3)
$$
giving rise to isomorphisms of $\O[G_K]$-modules

$$
D(A)(1) = T^*(1) = \Hom_{\zp}(T_{\wp}(B), \zp)(1) \iso T_{\wp}(\Bd), \qquad
A^*(1) = D(T)(1) \iso \Bd[\wp^\infty].
$$
Once we fix an isomorphism (1.3.1) via (1.3.2), we can pass from the Weil
pairing (1.3.3) to its $\O$-bilinear version, namely

$$
(\; ,\;)_{\O} : T_{\wp}(B) \times T_{\wp}(\Bd) = T \times T^*(1) \lo \O(1),
\qquad (x, y) = \Tr_{\K/\qp} (a (x, y)_{\O}),
$$
which induces an isomorphism of $\O[G_K]$-modules

$$
T^*(1) = \Hom_{\O}(T_{\wp}(B), \O)(1) \iso T_{\wp}(\Bd).
$$
The symmetric isogeny $\lambda$ from (1.1) defines morphisms of
$\O[G_K]$-modules

$$
\lambda_* : T_{\wp}(B) \ho T_{\wp}(\Bd),\qquad
B[\wp^\infty] \twoheadrightarrow \Bd[\wp^\infty]
$$
with finite cokernel and kernel, respectively. The Weil pairing attached
to $\lambda$

$$
(\; ,\;)_{\O,\lambda} : T_{\wp}(B) \times T_{\wp}(B) = T \times T \lo \O(1),
\qquad (x, y)_{\O,\lambda} := (x, \lambda_*(y))_{\O}
\eqno (1.3.4)
$$
is skew-symmetric; in other words, $\lambda_* : T \lo T^*(1)$ satisfies
$(\lambda_*)^*(1) = -\lambda_*$.

\sn
{\bf (1.4) Classical Selmer groups.} \ For every finite extension $L/K$,
$p$-power descent on $B$ over $L$ gives rise to the classical Selmer groups
$\Sel_{p^k}(B/L) \subset \coh 1L{B[p^k]}$ sitting in the standard
exact sequences

$$
0 \lo B(L) \otimes \Z/p^k\Z \lo \Sel_{p^k}(B/L) \lo \sha(B/L)[p^k] \lo 0.
\eqno (1.4.1)
$$
Their respective inductive and projective limits

$$
\Sel_{p^\infty}(B/L) := \dirlim_k \Sel_{p^k}(B/L) \subset \coh 1L{B[p^\infty]},
\qquad S_p(B/L) := \invlim_k \Sel_{p^k}(B/L) \subset \coh 1L{T_p(B)}
$$
coincide with the corresponding Bloch--Kato Selmer groups

$$
\cohbk L{B[p^\infty]} \subset \coh 1L{B[p^\infty]},\qquad
\cohbk L{T_p(B)} \subset \coh 1L{T_p(B)},
$$
by [BK, (3.11.1), (3.11.2)]. All groups in (1.4.1) and
in the limit exact sequences

$$
\displaylines{
0 \lo B(L) \otimes \qp/\zp \lo \Sel_{p^\infty}(B/L) \lo \sha(B/L)[p^\infty]
\lo 0,\cr
0 \lo B(L) \otimes \zp \lo S_p(B/L) \lo T_p \sha(B/L)[p^\infty] \lo 0\cr}
$$
are $O_M \otimes \zp$-modules. After tensoring with $\O$ over
$O_M \otimes \zp$, we obtain exact sequences

$$
0 \lo B(L) \otimes_{O_M} O_M/\wp^{ke} \lo \Sel_{\wp^{ke}}(B/L) \lo
\sha(B/L)[\wp^{ke}] \lo 0
\eqno (1.4.2)
$$
(where $e = e_{\wp}$ is the ramification index of $\wp$ above $p$) and

$$
\displaylines{
0 \lo B(L) \otimes_{O_M} \K/\O \lo \Sel_{\wp^\infty}(B/L) \lo
\sha(B/L)[\wp^\infty] \lo 0,\cr
0 \lo B(L) \otimes_{O_M} \O \lo S_{\wp}(B/L) \lo T_{\wp} \sha(B/L)[p^\infty]
\lo 0.\cr}
$$
Again,

$$
\Sel_{\wp^\infty}(B/L) = \cohbk L{B[\wp^\infty]} = \cohbk LA,\qquad
S_{\wp}(B/L) = \cohbk L{T_{\wp}(B)} = \cohbk LT.
$$
The same discussion applies to $\Bd$; one obtains

$$
\Sel_{\wp^\infty}(\Bd/L) = \cohbk L{\Bd[\wp^\infty]} = \cohbk L{A^*(1)},
\qquad S_{\wp}(\Bd/L) = \cohbk L{T_{\wp}(\Bd)} = \cohbk L{T^*(1)}.
$$
If $L \subset \o K$ is an arbitrary algebraic extension of $K$, we let

$$
\Sel_{\wp^{ke}}(B/L) := \dirlim_{F, \res} \Sel_{\wp^{ke}}(B/F)
\qquad (k \in \N \cup \{\infty\}),\qquad
S_{\wp}(B/L) := \invlim_{F, \cor} S_{\wp}(B/F),
$$
where $F$ runs through all intermediate fields $K \subset F \subset L$ such
that $[F : K] < \infty$.

\sn
{\bf (1.5) Greenberg's Selmer groups.} \ Let $v \mid p$ be a prime of $K$
above $p$. As $B$ has good $\wp$-ordinary reduction at $v$, there are
exact sequences of $\O[G_{K_v}]$-modules

$$
0 \lo T_v^+ \lo T \lo T_v^- \lo 0,\qquad
0 \lo A_v^+ \lo A \lo A_v^- \lo 0
\eqno (1.5.1)
$$
in which

$$
T_v^- = T_\wp(\widetilde B_v),\qquad A_v^- = \widetilde B_v[\wp^\infty],
$$
and the Pontryagin dual of (1.5.1) is isomorphic to

$$
0 \lo A^*(1)_v^+ \lo A^*(1) \lo A^*(1)_v^- \lo 0,\quad
0 \lo T^*(1)_v^+ \lo T^*(1) \lo T^*(1)_v^- \lo 0,
$$
where

$$
T^*(1)_v^- = T_\wp(\widetilde B_v^t),\qquad
A^*(1)_v^- = \widetilde B_v^t[\wp^\infty].
$$
In addition, $\lambda : B \lo \Bd$ induces maps

$$
T \ho T^*(1),\qquad T_v^\pm \ho T^*(1)_v^\pm,\qquad
A \twoheadrightarrow A^*(1),\qquad A_v^\pm \twoheadrightarrow A^*(1)_v^\pm
$$
with finite cokernel (for $T, T_v^\pm$) and kernel (for $A, A_v^\pm$),
respectively.

Fix a finite set $S$ of primes of $K$ containing all archimedean primes,
all primes above $p$ and all primes at which $B$ has bad reduction.
If $L$ is a finite extension of $K$, let $L_S$ be the maximal algebraic
extension of $L$ unramified outside primes above $S$; set $G_{L,S} :=
\Gal(L_S/L)$. Denote by $\Sigma_L$ (resp. $\Sigma^\prime_L$) the set
of all primes of $L$ above $p$ (resp. the set of all nonarchimedean primes
of $L$ above $S \smallsetminus \Sigma_K$). For each $X = T, A, T^*(1), T^*(1)$,
the Greenberg Selmer group over $L$ and its strict counterpart are defined,
respectively, by

$$
\eqalign{
S_X(L) &:= \Ker\(\coh 1{G_{L,S}}X \lo \bigoplus_{v\in\Sigma_L}
\coh 1{I_v}{X_v^-} \oplus \bigoplus_{v\in\Sigma_L^\prime} \coh 1{I_v}X\)\cr
S^\str_X(L) &:= \Ker\(\coh 1{G_{L,S}}X \lo \bigoplus_{v\in\Sigma}
\coh 1{G_{L_v}}{X_v^-} \oplus \bigoplus_{v\in\Sigma_L^\prime}
\coh 1{I_v}X\),\cr}
$$
where $I_v \subset G_{L_v} = \Gal(\o L_v/L_v)$ denotes the inertia group
at $v$. These groups do not depend on $S$, and the morphisms

$$
S_T(L) \otimes_{\O} \K/\O \ho S_A(L),\qquad
S_{T^*(1)}(L) \otimes_{\O} \K/\O \ho S_{A^*(1)}(L)
$$
(as well as their strict counterparts) have finite cokernels.

\sn
{\bf (1.6) Selmer complexes and extended Selmer groups.} \ 
In the notation of 1.5, the Selmer complex attached to $X = T, A, T^*(1),
T^*(1)$ over $L$ is defined as

$$
\Cff LX = \Cone\(C^\ddt_{\cont}(G_{L,S}, X) \oplus \bigoplus_{v\in \Sigma_L
\cup \Sigma_L^\prime} U_v^+(X) \lo \bigoplus_{v\in \Sigma_L \cup
\Sigma_L^\prime} C^\ddt_{\cont}(G_{L_v}, X)\) [-1],
$$
where

$$
U_v(X)^+ = \cases{
C^\ddt_{\cont}(G_{L_v}, X_v^+), & \qquad $v\in \Sigma_L$\cr
C^\ddt_{\cont}(G_{L_v}/I_v, X^{I_v}), & \qquad $v\in \Sigma^\prime_L$.\cr}
$$
Up to a canonical quasi-isomorphism, $\Cff LX$ does not depend on $S$;
its cohomology groups are denoted by $\cohf iLX$.

\sn
{\bf (1.7) Comparison of Selmer groups.} \ For each $X = T, A, T^*(1), T^*(1)$,
there is an exact sequence

$$
0 \lo \cohf 0LX \lo \coh 0LX \lo \bigoplus_{v\in \Sigma_L} \coh 0{L_v}{X_v^-}
\lo \cohf 1LX \lo S^\str_X(L) \lo 0,
$$
by [N2, Lemma 9.6.3]. In addition, [N2, Lemma 9.6.7.3] implies that there are
exact sequences

$$
0 \lo S^\str_T(L) \lo S_{\wp}(B/L) \lo \bigoplus_{v\in\Sigma_L}
{\coh 1{L_v}{T^-_v}}_\tors \oplus \bigoplus_{v\in\Sigma_L^\prime}
\coh 1{L_v}{T}/\cohur {L_v}{T}
$$

$$
\displaylines{
0 \lo \Sel_{\wp^\infty}(B/L) \lo S^\str_A(L) \lo\cr
\lo \bigoplus_{v\in\Sigma_L} \Imm\(\coh 1{L_v}{A^+_v} \lo
\coh 1{L_v}{A}\)/\divv \oplus
\bigoplus_{v\in\Sigma_L^\prime} \cohur {L_v}{A},\cr}
$$
in which

$$
\displaylines{
{\coh 1{L_v}{T^-_v}}_\tors \iso \coh 0{L_v}{A^-_v}/\divv =
\widetilde B_v(k(v))[\wp^\infty]\cr
D\(\Imm\(\coh 1{L_v}{A^+_v} \lo \coh 1{L_v}{A}\)/\divv\) \subseteq
\coh 0{L_v}{A^*(1)^-_v}/\divv = \widetilde B_v^t(k(v))[\wp^\infty].\cr}
$$
The same lemma implies that, for each $v\in\Sigma_L^\prime$, the $\O$-modules
$\coh 1{L_v}{T}/\cohur {L_v}{T}$ and $\cohur {L_v}{A}$ have the same finite
length, equal to the local Tamagawa factor $\Tam vT\wp$ defined in 1.8 below.

Of course, one can replace $B$ by $\Bd$, $T$ by $T^*(1)$ and $A$ by $A^*(1)$
in the above discussion.

\sn
{\bf (1.8) Local Tamagawa factors.} \ In the notation of 1.6, if $v \nmid p$
is a finite prime of $L$, the local Tamagawa factor $\Tam vT\wp$ is defined
as in [N2, 7.6.10] (following [FoPR, Prop. 4.2.2(ii)]), namely

$$
\Tam vT\wp := \ell_{\O}\(H^1(I_v,T)_\tors^{Fr(v) = 1}\)
$$
(where $Fr(v)$ is the arithmetic Frobenius at $v$ and $\ell_{\O}(Z)$ denotes
the length of any $\O$-module $Z$). This
is a non-negative integer (since the group $H^1(I_v,T)_\tors \simeq
H^0(I_v,A)/\divv$ is finite), equal to zero if $v\not\in \Sigma_L^\prime$.

It will be more convenient to use geometric notation; let us write

$$
\Tam v{B/L}\wp := \Tam vT\wp,\qquad
{\rm Tam}(B/L,\wp) := \sum_{v\in\Sigma_L^\prime} \Tam v{B/L}\wp.
$$
The equality

$$
\Tam v{T^*(1)}\wp = \Tam vT\wp
\eqno (1.8.1)
$$
proved in [N2, 10.2.8] then implies that

$$
\Tam v{\Bd/L}\wp := \Tam v{T^*(1)}\wp = \Tam v{B/L}\wp,\qquad
{\rm Tam}(\Bd/L,\wp) = {\rm Tam}(B/L,\wp).
$$
This cohomological definition agrees with the geometric one, namely, that

$$
\Tam v{B/L}\wp := \ell_{\O}\(\pi_0(\widetilde B_v)^{G_{k(v)}}
\otimes_{O_M}\O\).
\eqno (1.8.2)
$$
In particular, if $M = \Q$, then $\wp = p$ and $\Tam v{B/L}\wp$ is
equal to the $p$-adic valuation of the usual local Tamagawa factor
$c_{{\rm Tam},v}(B/L) = \# \coh 0{k(v)}{\pi_0(\widetilde B_v)}$.

Note that (1.8.2) also implies (1.8.1), by the $G_{k(v)}$-equivariance and
nondegeneracy of Grothendieck's monodromy pairing $\pi_0(\widetilde B_v)
\times \pi_0(\widetilde B_v^t) \lo \Q/\Z$.

\vskip16pt

\centerline{\bf 2. Comparison of Selmer groups, duality, control theorems}

\bn
{\bf (2.1) Conditions on $B$.} \ Given a finite extension $L/K$, consider
the following conditions.

\medskip
\item{$(A1)_{B,L,\wp}$} \quad There is an isomorphism of $\O[G_L]$-modules
$j : T\iso T^*(1)$ (where $T = T_{\wp}(B)$) such that $j^*(1) = -j$.
\item{$(A2)_{B,L,\wp}$} \quad ${\rm Tam}(B/L,\wp) = 0$ and
$\bigoplus_{v\in\Sigma_L} \widetilde B_v(k(v))[\wp] = 0$.
\item{$(A3)_{B,L,\wp}$} \quad $B(L)[\wp] = 0$.
\item{$(A4)_{B,L,\wp}$} \quad $\Sel_{\wp^\infty}(B/L) \iso \K/\O$ (which is
equivalent to $D(\Sel_{\wp^\infty}(B/L)) \iso \O$).

\sn

\proclaim{(2.2) Proposition}. (1) The conditions $(A2)_{B,L,\wp}$ and
$(A2)_{\Bd,L,\wp}$ are equivalent.\hb
(2) If $k \in \{3,4\}$ and if the conditions $(A1)_{B,L,\wp}$ and
$(Ak)_{B,L,\wp}$ hold, so does $(Ak)_{\Bd,L,\wp}$.\hb
(3) If $(A3)_{B,L,\wp}$ holds, then $\Sel_{\wp^{ke}}(B/L) =
\Sel_{\wp^\infty}(B/L)[\wp^{ke}]$ holds, for all $k\geq 1$.\hb
(4) If $L^\prime/L$ is a finite extension of $p$-power degree which is
unramified at all primes of $L$ at which $B$ has bad reduction, then
the conditions $(A2)_{B,L,\wp}$ and $(A2)_{B,L^\prime,\wp}$ are equivalent.\hb
(5) If $L^\prime/L$ is a finite Galois extension of $p$-power degree, then
the conditions $(A3)_{B,L,\wp}$ and $(A3)_{B,L^\prime,\wp}$ are equivalent.\hb
(6) If $\dim(B) = [M:\Q]$, then $(A1)_{B,L,\wp}$ holds, and the isomorphism
$ j : T\iso T^*(1)$ induces isomorphisms of $\O[G_{L_v}]$-modules $X_v^\pm \iso
X^*(1)_v^\pm$, for $X = T, A$ and all $v\in \Sigma_L$.\hb

\sn
{\it Proof.\/} (1) Combine (1.8.1) with the fact that
$\widetilde B_v(k(v))[\wp^\infty]$ and $\widetilde B_v^t(k(v))[\wp^\infty]$
have the same cardinality. The statement (2) is immediate, while (3) follows
from (1.4.2).
The statements (4) and (5) are consequences of the fact that, if a $p$-group
$G$ acts on a finite set $X$, then $\#(X^G) \equiv \#(X) \; ({\rm mod}\ p)$.

In the situation of (6), the given $O_M$-linear symmetric isogeny
$\lambda = \lambda^t : B \lo \Bd$ induces an isomorphism of $\K[G_K]$-modules
$\lambda_* : T_{\wp}(B)\otimes_{\O} \K \iso T_{\wp}(\Bd)\otimes_{\O} \K$ and
a $G_K$-equivariant, $\O$-bilinear, skew-symmetric pairing
$\langle \;,\;\rangle_{\O,\lambda}$ from (1.3.4), which is
non-degenerate when tensored with $\K$ and which satisfies
$T_{\wp}(\Bd) = \lambda_* \{y \in T_{\wp}(B) \otimes_{\O} \K \mid \forall
x\in T_{\wp}(B)\;\; \langle x,y\rangle_{\O,\lambda} \in \O\}$.

As $T = T_{\wp}(B)$ is a free $\O$-module of rank two, the matrix of the pairing
$\langle \;,\;\rangle_{\O,\lambda}$ in any basis of $T$ over $\O$ is of the form
$\pmatrix{0 & b\cr -b & 0\cr}$, for some $b\in \O\smallsetminus \{0\}$. 
This implies that $T_{\wp}(\Bd) = b^{-1} \lambda_*(T_{\wp}(B))$, hence
$j := b^{-1} \circ \lambda_* : T_{\wp}(B) \iso T_{\wp}(\Bd)$ has the required
property.

The maps $X^\pm_v \lo X^*(1)^\pm_v$ are isomorphisms, since $T_v^-$ is the
unique quotient of $T$ which is free of rank one over $\O$ on which $I_v$
acts trivially.

\proclaim{(2.3) Proposition}. Let $L$ be a finite extension of $K$.\hb
(1) If $(A2)_{B,L,\wp}$ holds, then
$$
\eqalign{
\Sel_{\wp^\infty}(B/L) &= \cohbk LA = S^\str_A(L) = \cohf 1LA,\qquad\qquad
\qquad\;\; A = B[\wp^\infty]\cr
\Sel_{\wp^\infty}(\Bd/L) &= \cohbk L{A^*(1)} = S^\str_{A^*(1)}(L) =
\cohf 1L{A^*(1)},\qquad A^*(1) = \Bd[\wp^\infty]\cr
S_\wp(B/L) &= \cohbk LT = S^\str_T(L) = \cohf 1LT,\qquad\qquad\qquad\;\;
T = T_\wp(B)\cr
S_\wp(\Bd/L) &= \cohbk L{T^*(1)} = S^\str_{T^*(1)}(L) = \cohf 1L{T^*(1)},\qquad
T^*(1) = T_\wp(\Bd).\cr}
$$
(2) In general (without assuming any $(Ak)_{B,L,\wp}$), there are isomorphisms
of $\O$-modules
$$
\eqalign{
D(\cohf iLT) &\simeq \cohf {3-i}L{A^*(1)} \qquad\qquad
(= 0 {\rm\ if\ } i \not= 1,2,3),\cr
D(\cohf iLA) &\simeq \cohf {3-i}L{T^*(1)} \qquad\qquad
(= 0 {\rm\ if\ } i \not= 0,1,2).\cr}
$$
(3) If $(A2)_{B,L,\wp}$ holds, then
$$
\cohf iLA = \cases{
{\rm a\ submodule\ of\ } B(L)[\wp^\infty], & $i=0$\cr
\Sel_{\wp^\infty}(B/L), & $i=1$\cr
D(S_\wp(\Bd/L)), & $i=2$\cr
0, & $i\not= 0,1,2$\cr}
$$
$$
\cohf iLT = \cases{
S_\wp(B/L), & $i=1$\cr
D(\Sel_{\wp^\infty}(\Bd/L)), & $i=2$\cr
{\rm a\ quotient\ of\ } D(\Bd(L)[\wp^\infty]), & $i=3$\cr
0, & $i\not= 1,2,3$\cr}
$$
(4) If $(A3)_{B,L,\wp}$ holds, then $\cohf 0LA = 0 = \cohf3L{T^*(1)}$.
Dually, if $(A3)_{\Bd,L,\wp}$ holds, then $\cohf 0L{A^*(1)} = 0 = \cohf3LT$.

\sn
{\it Proof.\/} The equalities of the various Selmer groups in (1) follow
from the discussion in 1.7. The statement (2) is a consequence of
[N2, Thm. 6.3.4, Prop. 6.7.7], while (3) is a combination of (1) and (2).
Finally, (4) follows from (2) and the fact that $\cohf 0LA \subset H^0(L,A)$.

\sn
{\bf (2.4) Iwasawa theory.} \ Fix a Galois extension $\ki/K$ such that
$\Gamma := \Gal(\ki/K) \simeq {\bf Z}_p^d$ ($d\geq 1$) and let $\Lambda :=
\O[[\Gamma]] = \invlim_F \O[\Gamma_F]$, where $F$ runs through all fields
$K \subset F \subset \ki$ such that $[F : K] < \infty$, and $\Gamma_F :=
\Gal(F/K)$.

For every intermediate field $K \subset L \subset \ki$ (not necessarily
of finite degree over $K$), let

$$
\Gamma^L := \Gal(\ki/L),\qquad
\Gamma_L := \Gal(L/K) = \Gamma/\Gamma^L,\qquad
\Lambda_L := \O[[\Gamma_L]].
$$
The corresponding Iwasawa-theoretical Selmer modules

$$
\cohf iLA := \dirlim_{F, \res} \cohf iFA,\qquad
\cohiwf iLKT := \invlim_{F, \cor} \cohf iFT
$$
($K \subset F \subset L$, $[F : K] < \infty$) are $\Lambda_L$-modules
of cofinite and finite type, respectively.

The standard involution $\iota : \Lambda_L \lo \Lambda_L$ is induced
by the inverse map $\Gamma_L \lo \Gamma_L$, $\gamma \mapsto \gamma^{-1}$.
For any $\Lambda_L$-module $N$ we denote by $N^\iota$ the $\Lambda_L$-module
equal to $N$ as an $\O$-module, but on which every $r \in \Lambda_L$
acts as $\iota(r)$ does on $N$. Note that $\iota$ induces an isomorphism
of $\Lambda_L$-modules $\iota : \Lambda_L \iso \Lambda_L^\iota$.

This involution appears naturally when one compares Pontryagin duality
between $\Lambda_L$-modules of finite type (compact) and cofinite
type (discrete), defined by

$$
D_{\Lambda_L}(N) := \Hom_{\cont,\zp}(N, \Q_p/\Z_p),\quad
(rf)(n) := f(rn)\quad (r\in \Lambda_L,\, f\in D_{\Lambda_L}(N),\, n\in N),
$$
with Pontryagin duality for $\O$-modules with a continuous linear
action of $\Gamma_L$: in this case

$$
D(N) := \Hom_{\cont,\zp}(N, \Q_p/\Z_p),\quad
(\gamma\cdot f)(n) := f(\gamma^{-1}(n))\quad
(\gamma\in \Gamma_L,\, f\in D(N),\, n\in N).
$$
In other words,

$$
D_{\Lambda_L}(N) = D(N)^\iota.
$$

\proclaim{(2.5) Proposition}. Assume that $K \subset L \subset \ki$ is
an arbitrary intermediate field.\hb
(1) In general (without assuming any $(Ak)_{B,L,\wp}$), there are isomorphisms
of $\Lambda_L$-modules
$$
\eqalign{
D_{\Lambda_L}(\cohiwf iLKT) &\simeq \cohf {3-i}L{A^*(1)}^\iota \qquad\qquad
\quad\; (= 0 {\rm\ if\ } i \not= 1,2,3),\cr
D_{\Lambda_L}(\cohf iLA) &\simeq \cohiwf {3-i}LK{T^*(1)}^\iota \qquad\qquad
(= 0 {\rm\ if\ } i \not= 0,1,2).\cr}
$$
(2) If $(A2)_{B,L,\wp}$ holds, then
$$
\cohf iLA = \cases{
{\rm a\ submodule\ of\ } B(L)[\wp^\infty], & $i=0$\cr
\Sel_{\wp^\infty}(B/L), & $i=1$\cr
D_{\Lambda_L}(S_\wp(\Bd/L))^\iota, & $i=2$\cr
0, & $i\not= 0,1,2$\cr}
$$
$$
\cohiwf iLKT = \cases{
S_\wp(B/L), & $i=1$\cr
D_{\Lambda_L}(\Sel_{\wp^\infty}(\Bd/L))^\iota, & $i=2$\cr
{\rm a\ quotient\ of\ } D_{\Lambda_L}(\Bd(L)[\wp^\infty])^\iota, & $i=3$\cr
0, & $i\not= 1,2,3$\cr}
$$
(3) (Exact control theorem) If $(A2)_{B,L,\wp}$ and $(A3)_{B,L,\wp}$ hold,
then the canonical map
$$
\Sel_{\wp^\infty}(B/L) \iso \Sel_{\wp^\infty}(B/\ki)^{\Gamma^L}
$$
is an isomorphism (idem if we replace everywhere $B$ by $\Bd$).\hb
(4) If $(A2)_{B,L,\wp}$ and $(A3)_{B,L,\wp}$ are satisfied, then there
is an exact sequence of $\Lambda_L$-modules of cofinite type
$$
\displaylines{
0 \lo H^1(\Gamma^L, \Sel_{\wp^\infty}(B/\ki)) \lo
D_{\Lambda_L}(S_\wp(\Bd/L))^\iota \lo
D_{\Lambda_L}(S_\wp(\Bd/\ki)_{\Gamma^L})^\iota\cr
\lo H^2(\Gamma^L, \Sel_{\wp^\infty}(B/\ki)) \lo 0\cr}
$$
(again, we can interchange everywhere $B$ with $\Bd$).\hb
(5) If $(A2)_{B,L,\wp}$ and $(A3)_{B,L,\wp}$ are satisfied, then there is an
isomorphism of $\Lambda_L$-modules of finite type
$$
S_\wp(B/L) \iso \Hom_{\Lambda_L}(D_{\Lambda_L}(\Sel_{\wp^\infty}(B/L)),
\Lambda_L)
$$
(as before, we can replace everywhere $B$ by $\Bd$).\hb
(6) If $(A2)_{B,L,\wp}$ and $(A3)_{B,L,\wp}$ are satisfied and if
$\Gamma_L \simeq {\bf Z}_p^r$ ($0 \leq r\leq d$), then
$$
\rk_{\Lambda_L}\, S_\wp(B/L) = \rk_{\Lambda_L}\, S_\wp(\Bd/L) =
\cork_{\Lambda_L}\, \Sel_{\wp^\infty}(B/L) =
\cork_{\Lambda_L}\, \Sel_{\wp^\infty}(\Bd/L).
$$

\sn
{\it Proof.\/} (1), (2) Apply Proposition 2.3(2)-(3) over all intermediate
fields $K \subset F \subset L$ such that $[F : K] < \infty$ (which is
legitimate, thanks to Proposition 2.2(4)-(5)) and take the inductive
(resp.  the projective) limit.

\noindent
(3),(4) In the spectral sequence from [N2, Prop. 8.10.12]

$$
{}^\prime\! \o E_2^{i,j} = H^i(\Gamma^L, \cohf j\ki{A}) \Lo \cohf {i+j}LA
$$
(which is a consequence of the ``exact control theorem for Selmer complexes"
[N2, Prop. 8.10.1]) we have ${}^\prime\! \o E_2^{i,j} = 0$ if $j\not= 1,2$,
by (1) applied to $\ki$ and the fact that $B(\ki)[\wp] = 0$ (which follows from
$(A3)_{B,K,\wp}$, by Proposition 2.2(5)).

\noindent
(5) The duality theorem [N2, Thm. 8.9.12] applies in this case, giving rise
to a spectral sequence

$$
E_2^{i,j} = \Ext^i_{\Lambda_L}(\cohiwf {3-j}LK{T^*(1)}, \Lambda_L)^\iota \Lo
\cohiwf {i+j}LKT
$$
satisfying $E_2^{i,j} = 0$ for $j\not= 1,2$ (as in the proof of (3) and (4)).
Therefore

$$
S_\wp(B/L) = \cohiwf 1LKT \simeq E_2^{0,1} \simeq
\Hom_{\Lambda_L}(D_{\Lambda_L}(\Sel_{\wp^\infty}(B/L)), \Lambda_L).
$$
(6) This follows from (5) and the fact that there exists a constant $C\geq 0$
such that, for every $k\geq 1$ and every finite extension $F/K$, the kernel
and the cokernel of the map

$$
\Sel_{\wp^k}(B/F) \lo \Sel_{\wp^k}(\Bd/F)
$$
induced by $\lambda : B \lo \Bd$ is killed by $\wp^C$.

\sn
{\bf (2.6) Notation.} \ For every field $K \subset L \subset \ki$ we are
going to abbreviate as

$$
0 \lo Z(B/L) \lo X(B/L) \lo Y(B/L) \lo 0
\eqno (2.6.1)
$$
the terms in the exact sequence

$$
0 \lo D_{\Lambda_L}(\sha(B/L)[\wp^\infty]) \lo
D_{\Lambda_L}(\Sel_{\wp^\infty}(B/L)) \lo
D_{\Lambda_L}(B(L) \otimes_{O_M} \K/\O) \lo 0.
$$
Proposition 2.5(3) tells us that, under the conditions $(A2)_{B,K,\wp}$ and
$(A3)_{B,K,\wp}$, there are canonical isomorphisms of $\Lambda_L$-modules

$$
X(B/\ki)_{\Gamma^L} \iso X(B/L),
\eqno (2.6.2)
$$
hence also

$$
X(B/L^\prime)_{\Gal(L^\prime/L)} \iso X(B/L)\qquad\qquad
(K \subset L \subset L^\prime \subset \ki).
$$

\vskip16pt

\centerline{\bf 3. Freeness of compact Selmer groups and the vanishing
of $\sha[\wp^\infty]$}

\bn
{\bf (3.1)} \ Consider another condition on $B$ and $K$.

\medskip

\item{$(A5)_{B,K}$} \quad There exists a subfield $K^+ \subset K$ such that
$[K:K^+] = 2$, $\ki/K^+$ is a Galois extension with Galois group $\Gamma^+ =
\Gamma \rtimes \{1, c\}$, where $c^2 = 1$ and $\forall \gamma\in\Gamma
\;\; c \gamma c^{-1} = \gamma^{-1}$, and there exists an abelian variety
$B^+$ over $K^+$ with good reduction at all primes of $K^+$ above $p$,
equipped with a ring morphism $i^+ : O_M \lo \End(B^+)$ and a symmetric
$O_M$-linear isogeny $\lambda^+ = (\lambda^+)^t : B^+ \lo (B^+)^t$, such that
$B^+$ has good $\wp$-ordinary reduction at all primes of $K^+$ above $p$,
and that the base change of $(B^+, i^+, \lambda^+)$ from $K^+$ to $K$ is
isomorphic to $(B, i, \lambda)$.

\sn

\proclaim{(3.2) Proposition}. (1) If $(A5)_{B,K}$ holds, then
$$
\displaylines{
\cork_{\O}\, \Sel_{\wp^\infty}(B/K) \equiv
\cork_{\Lambda_L}\, \Sel_{\wp^\infty}(B/\ki) \; ({\rm mod}\ 2),\cr
\cork_{\O}\, \Sel_{\wp^\infty}(B/K) \geq
\cork_{\Lambda_L}\, \Sel_{\wp^\infty}(B/\ki).\cr}
$$
(2) If the conditions $(A2)_{B,K,\wp}$, $(A3)_{B,K,\wp}$ and $(A4)_{B,K,\wp}$
are satisfied, then, for every intermediate field $K \subset L \subset \ki$
such that $\Gamma_L = \Gal(L/K) \simeq {\bf Z}_p^r$ ($0\leq r\leq d$),
$X(B/L) = D_{\Lambda_L}(\Sel_{\wp^\infty}(B/L))$ is a cyclic
$\Lambda_L$-module.\hb
(3) If the conditions $(A2)_{B,K,\wp}$, $(A3)_{B,K,\wp}$, $(A4)_{B,K,\wp}$ and
$(A5)_{B,K}$ are satisfied, then, for every intermediate field $K \subset L
\subset \ki$, the $\Lambda_L$-modules $X(B/L)$ and $S_\wp(B/L)$ are free
of rank one, the $\Lambda_L/\wp^{ke}$-module
$D_{\Lambda_L}(\Sel_{\wp^{ke}}(B/L))$ is free of rank one
(for every $k\geq 1$), and the canonical maps
$$
X(B/L^\prime)_{\Gal(L^\prime/L)} \iso X(B/L),\qquad
S_\wp(\Bd/L^\prime)_{\Gal(L^\prime/L)} \iso S_\wp(\Bd/L)
$$
are isomorphisms of $\Lambda_L$-modules, if $K \subset L \subset
L^\prime \subset \ki$. In particular, if $[L : K] < \infty$,
then
$$
\rk_{O_M}\, B(L) + \cork_{\O}\, \sha(B/L)[\wp^\infty] = [L : K].
$$

\sn
{\it Proof.\/} (1) If $B = E$ is an elliptic curve, this is
[N2, Prop. 10.7.19]. The general case follows from [N2, Thm. 10.7.17(iv)]
(for $\chi = \chi^\prime = 1$).

\noindent
(2) $X(B/L)$ is a $\Lambda_L$-module of finite type satisfying
$X(B/L)_{\Gamma_L} \iso X(B/K) \iso \O$, by (2.6.2) (for $L/K$ replacing
$\ki/L$) and $(A4)_{B,K,\wp}$, respectively. If the image of $x\in X(B/L)$
generates $X(B/K)$ as an $\O$-module, then $(X(B/L)/\Lambda_L x)_{\Gamma_L}
= 0$, hence $X(B/L) = \Lambda_L x$ by Nakayama's Lemma.

\noindent
(3) It is enough to treat the case $L^\prime = \ki$.
According to (2) applied to $L = \ki$, we have $X(B/\ki) \iso \Lambda/J$
for some ideal $J \subset \Lambda$. On the other hand, (1) together with
$(A4)_{B,K,\wp}$ imply that $\rk_\Lambda\, X(B/\ki) = 1$, hence $J = 0$
and $X(B/\ki)$ is free of rank one over $\Lambda$. The control theorem
(2.6.2) then yields $X(B/L) \iso \Lambda_L$ as a $\Lambda_L$-module.
The statement about $\Sel_{\wp^{ke}}(B/L)$ then follows from Proposition
2.2(3).

It remains to show that $S_\wp(\Bd/\ki)_{\Gamma^L} \iso S_\wp(\Bd/L)$
is an isomorphism, which is equivalent, by Proposition 2.5(4), to the
vanishing of $H^i(\Gamma^L, \Sel_{\wp^\infty}(B/\ki))$ for $i = 1,2$.
We claim that the latter group vanishes for all $i > 0$. Indeed, its
Pontryagin dual $H_i(\Gamma^L, X(B/\ki)) \simeq H_i(\Gamma^L, \Lambda)$
is the $i$-th homology group of the Koszul complex of $\Lambda$ with
respect to the sequence $(\gamma_1^\prime - 1,\ldots, \gamma_t^\prime - 1)$,
where $\gamma_1^\prime, \ldots,\gamma_t^\prime$ is any basis of $\Gamma^L
\simeq {\bf Z}_p^t$ over $\zp$. We can take $\gamma_i^\prime =
\gamma_i^{p^{n_i}}$ ($1\leq i \leq t$, $n_i\geq 0$), for a suitable basis
$\gamma_1, \ldots,\gamma_d$ of $\Gamma$ over $\zp$. Therefore $\Lambda =
\O[[X_1,\ldots,X_d]]$ ($X_i = \gamma_i - 1$) and $\gamma_i^\prime - 1 =
(1 + X_i)^{p^{n_i}} - 1$, which implies that $(\gamma_1^\prime - 1,\ldots,
\gamma_t^\prime - 1)$ is a regular sequence in $\Lambda$, hence
$H_i(\Gamma^L, \Lambda) = 0$ for $i > 0$.

\sn
{\bf (3.3) Notation.} \ For an arbitrary algebraic extension $K^\prime/K$, let
us write

$$
N_{K^\prime/K}(B \otimes \O) := \invlim_{F, \Tr} (B(F) \otimes_{O_M} \O)
\subset S_\wp(B/K^\prime),
$$
where $F$ runs through all intermediate fields $K \subset F \subset K^\prime$
such that $[F : K] < \infty$.

We are going to consider the following conditions.

\medskip
\item{$(A6)_{B,K^\prime/K,\wp}$} \quad $\Imm(N_{K^\prime/K}(B \otimes \O) \lo
B(K) \otimes_{O_M} \O) \not= 0$.
\item{$(A7)_{B,K^\prime/K,\wp}$} \quad $\Imm(N_{K^\prime/K}(B \otimes \O) \lo
B(K) \otimes_{O_M} \O/\wp) \not= 0$.

\sn

\proclaim{(3.4) Theorem}. Assume that the conditions $(A2)_{B,K,\wp}$,
$(A3)_{B,K,\wp}$, $(A4)_{B,K,\wp}$, $(A5)_{B,K}$ and $(A6)_{B,\ki/K,\wp}$
are satisfied. Then, for every intermediate field $K \subset L \subset \ki$
such that $\Gamma_L = \Gal(L/K) \simeq {\bf Z}_p^r$ ($0\leq r\leq d$),
$$
\sha(B/L)[\wp^\infty] = 0,\qquad
B(L) \otimes_{O_M} \K/\O = \Sel_{\wp^\infty}(B/L)
$$
(the Pontryagin dual of the latter group being a free module of rank one
over $\Lambda_L$).

\sn
{\it Proof.\/} Induction on $r$.
If $r = 0$, then $L = K$. In this case $B(K) \otimes_{O_M} \O
\not= 0$ and $B(K)[\wp] = 0$ by $(A6)_{B,\ki/K,\wp}$ and $(A3)_{B,K,\wp}$,
respectively. This means that $B(K) \otimes_{O_M} \O \simeq \O^m$ and
$B(K) \otimes_{O_M} \K/\O \simeq (\K/\O)^m$ for some $m\geq 1$. On the other
hand, $B(K) \otimes_{O_M} \K/\O \subset \Sel_{\wp^\infty}(B/K) \simeq
\K/\O$ (by $(A4)_{B,K,\wp}$); thus $m = 1$ and $\sha(B/K)[\wp^\infty] = 0$.

Assume that $r > 0$. In the notation of 2.6, we need to show that $Z(B/L) = 0$,
which is equivalent to $X(B/L)/Z(B/L) = Y(B/L) = D_{\Lambda_L}(B(L)
\otimes_{O_M} \K/\O)$ not being $\Lambda_L$-torsion, since $X(B/L) \simeq
\Lambda_L$, by Proposition 3.2(3). Note that the canonical map

$$
B(F^\prime) \otimes_{O_M} \K/\O \lo
\left(B(F) \otimes_{O_M} \K/\O\right)^{\Gal(F/F^\prime)}
$$
is injective, whenever $K \subset F^\prime \subset F\subset \ki$ and
$[F : K] < \infty$, since $B(\ki)[\wp] = 0$ (by Proposition 2.2(5)
and $(A3)_{B,K,\wp}$).

If $r = 1$, write $L = \bigcup_{n\geq 1} K_n$, where $\Gal(K_n/K) \simeq
\Z/p^n\Z$. If $Y(B/L)$ were $\Lambda_L$-torsion, it would be a free
$\O$-module of finite type (since $Y(B/L)[\wp] = 0$), hence $B(L)
\otimes_{O_M} \O = B(K_m) \otimes_{O_M} \O$ for some $m\geq 1$.
This would imply that

$$
\forall k\geq 0\;\; \forall n\geq m\quad
N_{K_{n+k}/K_n}(B(K_{n+k}) \otimes_{O_M} \O) \subset
p^k (B(K_n) \otimes_{O_M} \O),
$$
hence $N_{L/K}(B \otimes \O) = 0$, which contradicts $(A6)_{B,\ki/K,\wp}$.

Assume that $r > 1$.
If $Y(B/L)$ were $\Lambda_L$-torsion, there would be $\gamma \in \Gamma_L
\smallsetminus \Gamma_L^p$ such that $(\gamma - 1) \not\in
{\rm Supp}_{\Lambda_L}(Y(B/L))$. The fixed field $L^\prime := L^{\gamma = 1}$
satisfies $\Gamma_{L^\prime} \simeq {\bf Z}_p^{r-1}$. By construction,
$Y(B/L)/(\gamma - 1)$ is a torsion module over $\Lambda_L/(\gamma - 1) =
\Lambda_{L^\prime}$, hence so is its quotient $Y(B/L^\prime)$; but this
is false by the induction hypothesis.

\proclaim{(3.5) Theorem}. Assume that the conditions $(A2)_{B,K,\wp}$,
$(A3)_{B,K,\wp}$, $(A4)_{B,K,\wp}$, $(A5)_{B,K}$ and $(A7)_{B,\ki/K,\wp}$
are satisfied.\hb
(a) For every intermediate field $K \subset L \subset \ki$ the following
statements hold.
$$
\sha(B/L)[\wp^\infty] = 0,\qquad
B(L) \otimes_{O_M} \K/\O = \Sel_{\wp^\infty}(B/L),\qquad
N_{L/K}(B \otimes \O) = S_\wp(B/L)
$$
and both $S_\wp(B/L)$ and $D_{\Lambda_L}(\Sel_{\wp^\infty}(B/L))$ are free
modules of rank one over $\Lambda_L = \O[[\Gal(L/K)]]$.
In the special case when $[L : K] < \infty$, then $B(L) \otimes_{O_M} \O =
S_\wp(B/L)$ is a free $\O$-module of rank $\rk_{O_M}\, B(L) = [L : K]$.\hb
(b) If, in addition, $(A1)_{B,K,\wp}$ is satisfied, then the canonical maps
$$
N_{L^\prime/K}(B \otimes \O) = S_\wp(B/L^\prime) \lo S_\wp(B/L) =
N_{L/K}(B \otimes \O)
$$
are surjective, for arbitrary intermediate fields $K \subset L \subset
L^\prime \subset \ki$. Furthermore, $S_\wp(B/L)$ is generated as a
$\Lambda_L$-module by the image $x_L$ of any element $x \in N_{\ki/K}
(B \otimes \O)$ whose image $\o x_K$ in $B(K) \otimes_{O_M} \O/\wp$
is non-zero.

\sn
{\it Proof.\/} (a) Fix intermediate fields $K \subset L \subset L^\prime
\subset \ki$. For every finite extension $F/K$,

$$
\Cone(\Cff FX \, \lomapr {\lambda_*}\, \Cff F{X^*(1)}) \qquad\qquad (X = T, A)
$$
is quasi-isomorphic to a complex of $\O/\wp^C$-modules, where
$\wp^C\,\Ker(\lambda)(\o K)[\wp^\infty] = 0$.
Therefore the kernel and cokernel of

$$
\cohf iFX \lo \cohf iF{X^*(1)} \qquad\qquad (X = T, A)
$$
is killed by $\wp^C$, for every $i$. Thus the same is true for the kernels
and cokernels of the maps

$$
S_\wp(B/L) \lo S_\wp(\Bd/L),\qquad \Sel_{\wp^\infty}(B/L) \lo 
\Sel_{\wp^\infty}(\Bd/L).
$$
Combined with the freeness results and the isomorphisms in Proposition 3.2(3),
this implies that the canonical map

$$
j_{L^\prime/L} : S_\wp(B/L^\prime)_{\Gal(L^\prime/L)} \lo S_\wp(B/L)
$$
is a morphism between two free $\Lambda_L$-modules of rank one, whose
kernel and cokernel is killed by $\wp^C$. Therefore $\Ker(j_{L^\prime/L}) = 0$
and the maps in the commutative diagram

$$
\normalbaselines
\matrix{
N_{L^\prime/K}(B \otimes \O)_{\Gal(L^\prime/L)} & \lomapr {k_{L^\prime/L}} & 
S_\wp(B/L^\prime)_{\Gal(L^\prime/L)}\cr
\mapd {i_{L^\prime/L}} && \mapd {j_{L^\prime/L}}\cr
N_{L/K}(B \otimes \O) & \lomapr {k_L} & S_\wp(B/L)\cr}
$$
satisfy

$$
\Ker(j_{L^\prime/L}) = 0 = \Ker(k_L),\qquad \wp^C\,\Coker(j_{L^\prime/L}) = 0.
$$
Moreover, $S_\wp(B/K) \iso \O$ (by $(A4)_{B,K,\wp}$) and
$\Coker(k_K \circ i_{L^\prime/K}) = 0$, by $(A7)_{B,\ki/K,\wp}$. As a result,
$j_{L^\prime/K}$ is an isomorphism and

$$
0 = \Coker(k_{L^\prime/K}) =
(S_\wp(B/L^\prime)/N_{L^\prime/K}(B \otimes \O))_{\Gal(L^\prime/K)},
$$
which implies that

$$
N_{L^\prime/K}(B \otimes \O) = S_\wp(B/L^\prime) \simeq \Lambda_{L^\prime}
\eqno (3.5.1)
$$
(for arbitrary $K \subset L^\prime \subset \ki$), by Nakayama's Lemma.

In the special case when $[L^\prime : K] < \infty$, it follows from (3.5.1)
that

$$
B(L^\prime) \otimes_{O_M} \O = S_\wp(B/L^\prime),\qquad
\rk_{O_M} B(L^\prime) = \rk_{\O}\, S_\wp(B/L^\prime) = [L^\prime : K].
$$
As a result, both $X(B/L^\prime)$ and $Y(B/L^\prime)$ in the exact sequence
(2.6.2) are free modules over $\O$ of the same rank $[L^\prime : K]$, hence
$Z(B/L^\prime) = 0$ and $\sha(B/L^\prime)[\wp^\infty] = 0$. This proves
the Theorem in the special case when $[L : K] < \infty$. The general case
follows by taking the inductive limit over all subfields of $L$ of finite
degree over $K$.

\noindent
(b) In this case the arguments in the proof of (a) go through if one replaces
the map $\lambda_*$ by the map $j_*$ induced by the isomorphism $j : T_\wp(B)
\iso T_\wp(B^t)$ from $(A1)_{B,K,\wp}$. The constant $C$ is then replaced by
zero, which means that the map 
$j_{L^\prime/L} : S_\wp(B/L^\prime)_{\Gal(L^\prime/L)} \lo S_\wp(B/L)$
is an isomorphism between two free $\Lambda_L$-modules of rank one.
If $\o x_K \not= 0$, then $S_\wp(B/\ki)/\Lambda x = 0$, by Nakayama's Lemma.
It follows that $S_\wp(B/\ki) = \Lambda x$ and, after applying $j_{\ki/L}$,
that $S_\wp(B/L) = \Lambda_L x_L$.

\vskip16pt

\centerline{\bf 4. An application to Heegner points}

\bn
{\bf (4.1) Ring class fields.} \ Let $K$ be an imaginary quadratic field
of discriminant $D_K$. Denote by $\eta_K : (\Z/|D_K|\Z)^\times \lo \{\pm 1\}$
the primitive quadratic character attached to $K$. For each prime $p\nmid
2 D_K$, we have $\eta_K(p) = \left({{D_K} \over p}\right)$; if $p\mid D_K$,
then $\eta_K(p) = 0$.

For any integer $m\geq 1$, denote by $O_m := \Z + m\O_K \subset O_K$ the order
of conductor $m$ in $K$ and by $H_m$ the ring class field of $K$ of conductor
$m$ ($H_1$ is the Hilbert class field of $K$).
The Galois groups of the intermediate extensions in the diagram

$$
\Q = K^+ \ho K \ho H_1 \ho H_m
$$
are as follows.

$$
G_m := \Gal(H_m/K) \simeq \Pic(O_m),\qquad
\Gal(H_m/\Q) = G_m \rtimes \{1,c\},\qquad
\forall g\in G_m\;\; c g c^{-1} = g^{-1}
$$
(where $c$ is complex conjugation) and there is an exact sequence

$$
\fr {\O_K^\times}{\Z^\times} \lo
\fr{(O_K \otimes \Z/m\Z)^\times}{(\Z/m\Z)^\times} \lo G_m \lo G_1 \lo 0.
$$
The first group in this sequence is cyclic of order

$$
u_K := \# (\O_K^\times/\Z^\times) = \cases{
3, & \quad $D_K = -3$\cr
2, & \quad $D_K = -4$\cr
1, & \quad $D_K \not= -3, -4$.\cr}
$$

\sn
{\bf (4.2) The anticyclotomic $\zp$-extension $\ki/K.$} \ For a fixed prime
number $p$, the tower of fields

$$
\Q = K^+ \ho K \ho H_1 \ho H_p \ho H_{p^2} \ho \cdots \ho
H_{p^\infty} := \bigcup_{n\geq 0} H_{p^n}
$$
has the following properties.

\medskip
\item{$\bullet$} \quad $\forall n\geq 1 \quad \Gal(H_{p^{n+1}}/H_{p^n})
\simeq \Z/p\Z$.
\item{$\bullet$} \quad If $p\not= 2$, then $\Gal(H_{p^\infty}/H_p) \simeq \zp$.
\item{$\bullet$} \quad $\Gal(H_1/K) \simeq \Pic(O_K) =  Cl_K$.
\item{$\bullet$} \quad $\Gal(H_p/H_1)$ is a cyclic group of order
$u_K^{-1}(p - \eta_K(p))$.
\item{$\bullet$} \quad The torsion subgroup
$\Delta := \Gal(H_{p^\infty}/K)_{\tors}$ is finite. Its fixed field
$\ki := (H_{p^\infty})^\Delta$ satisfies $\Gal(\ki/K) \simeq \zp$ and
$\Gal(\ki/\Q) = \Gal(\ki/K) \rtimes \{1,c\}$, as in (A5). Write
$\ki = \bigcup_{n\geq 0} K_n$, where $\Gal(K_n/K) \simeq \Z/p^n\Z$.

\medskip

\sn
{\bf (4.3) Heegner points.} \ Assume that

\medskip
\item{$\bullet$} \quad $E$ is an elliptic curve over $\Q$ of conductor $N$.
\item{$\bullet$} \quad $\varphi : X_0(N) \lo E$ is a modular parameterisation
of $E$ (sending $i\infty$ to the origin) of the smallest degree.
\item{$\bullet$} \quad $K$ is an imaginary quadratic field satisfying the
Heegner condition
$$
{\rm\ all\ primes\ dividing\ } N {\rm\ split\ in\ } K/\Q.
\leqno ({\rm Heeg})
$$

\medskip
\noindent
Fix an ideal ${\cal N} \subset O_K$ such that $O_K/{\cal N} \simeq \Z/N\Z$.
If $m\geq 1$ is an integer such that $(m, N) = 1$, then ${\cal N}_m :=
{\cal N} \cap O_m$ is an invertible ideal of $O_m$ satisfying
$O_m/{\cal N}_m \simeq {\cal N}_m^{-1}/O_m \simeq \Z/N\Z$.

The Heegner points of conductor $m$ on $X_0(N)$ and $E$, respectively, are
defined as

$$
x_m := [\C/O_m \lo \C/{\cal N}_m^{-1}] \in X_0(N)(H_m),\qquad
y_m := \varphi(x_m) \in E(H_m)
$$
(up to a sign, $y_m$ does not depend on the choice of ${\cal N}$).
The basic Heegner point on $E$ is defined as

$$
y_K := \Tr_{H_1/K} (y_1) \in E(K).
$$
A general modular parameterisation $\varphi^\prime : X_0(N) \lo E$ of $E$
(sending $i\infty$ to the origin) is obtained by composing $\varphi$ with
a non-trivial element $a \in \End(E) = \Z$. The Heegner points $y_m^\prime
:= \varphi^\prime(x_m)$ corresponding to $\varphi^\prime$ are therefore equal
to $y_m^\prime = a y_m$.

\sn
{\bf (4.4) Norm relations.} \ Fix a prime number $p \nmid N$ and let
$a_p := p + 1 - \# \widetilde E_p(\F_p)$. For any integer $m\geq 1$
relatively prime to $pN$, the Heegner points of conductors $m p^n$
on $E$ are related as follows [PR, 3.1, Prop. 1].

$$
\eqalign{
\forall n\geq 1 \qquad \Tr_{H_{mp^{n+1}}/H_{mp^n}} (y_{mp^{n+1}}) &=
a_p y_{mp^n} - y_{mp^{n-1}},\cr
\exists\, \sigma \in \Gal(H_m/K) \qquad
u_{K,m} \cdot \Tr_{H_{mp}/H_m} (y_{mp}) &= \cases{
a_p y_m, & \quad if $\eta_K(p) = -1$\cr
(a_p - \sigma) y_m, & \quad if $\eta_K(p) = 0$\cr
(a_p - \sigma - \sigma^{-1}) y_m, & \quad if $\eta_K(p) = 1$,\cr}
\cr
u_{K,m} &= \cases{
u_K, & \quad if $m = 1$,\cr
1, & \quad if $m > 1$.\cr}\cr
u_K\cdot \Tr_{H_p/K} (y_p) &= (a_p - 1 - \eta_K(p)) y_K.\cr}
$$

\sn
{\bf (4.5) Universal norms in the $p$-ordinary case.} \ Assume that $E$ has
good ordinary reduction at a prime number $p$ (which is equivalent to
$p \nmid N\cdot a_p$). In this case the polynomial defining the Euler factor
of $E$ at $p$ factors in $\zp[X]$ as

$$
\eqalignno{
X^2 - a_p X + p &= (X - \alpha_p)(X - \beta_p),\quad
\alpha_p \in {\bf Z}_p^\times,\quad \beta_p\in p{\bf Z}_p^\times,
& (4.5.1)\cr
&\alpha_p + \beta_p = a_p,\quad \alpha_p \beta_p = p.
& (4.5.2)\cr}
$$
%$$
%X^2 - a_p X + p = (X - \alpha_p)(X - \beta_p),\quad
%\alpha_p \in {\bf Z}_p^\times,\quad \beta_p\in p{\bf Z}_p^\times,\quad
%\alpha_p + \beta_p = a_p,\quad \alpha_p \beta_p = p.
%\eqno (4.5.1)
%$$
In addition, $|\iota(\alpha_p)| = |\iota(\beta_p)| = \sqrt p$, for every
embedding $\iota : \Q(\alpha_p) \ho \C$.

Define, for every integer $n\geq 0$,

$$
z_n := \alpha_p^{-n} y_{p^{n+1}} - \alpha_p^{-n-1} y_{p^n} \in
E(H_{p^{n+1}}) \otimes \zp.
$$
These elements are norm compatible, namely

$$
\forall n\geq 1 \quad \Tr_{H_{p^{n+1}}/H_{p^n}} (z_n) = z_{n-1}.
\eqno (4.5.3)
$$
In addition, the bottom element $z_0 = y_p - \alpha_p^{-1} y_0$ satisfies

$$
u_K\cdot \Tr_{H_p/K} (z_0) = (a_p - 1 - \eta_K(p)) y_K -
\alpha_p^{-1} (p - \eta_K(p)) y_K =
(\alpha_p - 1)(1 - \alpha_p^{-1} \eta_K(p)) y_K.
\eqno (4.5.4)
$$

\proclaim{(4.6) Proposition}. Assume that $p \nmid N\cdot a_p$.\hb
(1) The element $(\alpha_p - 1)(1 - \alpha_p^{-1} \eta_K(p))\, (y_K\otimes 1)
\in E(K) \otimes \zp$ is contained in
$$
u_K\cdot \Imm(N_{H_{p^\infty}/K}(E \otimes \zp) \lo E(K) \otimes \zp) \subset
u_K\cdot \Imm(N_{\ki/K}(E \otimes \zp) \lo E(K) \otimes \zp).
$$
(2) If $v$ runs through all primes of $K$ above $p$, then
$$
\prod_{v\mid p} \# \widetilde E_p(k(v)) \equiv
(1 - \alpha_p)(1 - \alpha_p \eta_K(p)) \; ({\rm mod}\ p).
$$
(3) If $a_p \not\equiv 1, \eta_K(p) \; ({\rm mod}\ p)$, then
$p \nmid \prod_{v\mid p} \# \widetilde E_p(k(v))$ and
$$
y_K \otimes 1 \in u_K\cdot \Imm(N_{H_{p^\infty}/K}(E \otimes \zp) \lo E(K)
\otimes \zp) \subset u_K\cdot \Imm(N_{\ki/K}(E \otimes \zp) \lo E(K) \otimes
\zp).
$$

\sn
{\it Proof.\/} (1) This is a consequence of the norm relations
(4.5.3) and (4.5.4).

\noindent
(2) The term on the left hand side is equal to
$(p + 1 - a_p) (p + 1 - \eta_K(p) a_p)$ if $\eta_K(p) \not= 0$, resp. to
$p + 1 - a_p$ if $\eta_K(p) = 0$. The claim follows from the fact that
$a_p \equiv \alpha_p \; ({\rm mod}\ p)$.

\noindent
(3) This is an immediate consequence of (1) and (2).

\sn
{\bf (4.7)} \ We are now ready to combine the abstract Iwasawa-theoretical
results of \S{1}-\S{3} with the norm relations summarised in Proposition 4.6.

\proclaim{(4.8) Theorem}. If $p \not= 2$ is a prime number such that\hb
(a) $E(K)[p] = 0$,\hb
(b) $p\nmid N \cdot a_p \cdot (a_p - 1) \cdot c_{\rm Tam}(E/\Q)$,\hb
(c) $y_K \not\in E(K)_{\tors}$,\hb
(d) $\rk_{\Z}\, E(K) = 1$ and $\sha(E/K)[p^\infty] = 0$,\hb
then $\sha(E/\ki)[p^\infty] = 0$ and the Pontryagin dual of $E(\ki) \otimes
\qp/\zp = \Sel_{p^\infty}(E/\ki)$ is a free module of rank one over
$\zp[[\Gal(\ki/K)]]$.

\proclaim{(4.9) Theorem}. If $p \not= 2$ is a prime number such that\hb
(a) $E(K)[p] = 0$,\hb
(b') $p\nmid N \cdot a_p \cdot (a_p - 1) \cdot (a_p - \eta_K(p))\cdot
c_{\rm Tam}(E/\Q)$,\hb
(c') $y_K \not\in p E(K)$,\hb
(d) $\rk_{\Z}\, E(K) = 1$ and $\sha(E/K)[p^\infty] = 0$,\hb
then, for every intermediate field $K \subset L \subset \ki$,
$\sha(E/L)[p^\infty] = 0$ and the Pontryagin dual of $E(L) \otimes \qp/\zp =
\Sel_{p^\infty}(E/L)$ is a free module of rank one over 
$\zp[[\Gal(L/K)]]$. For every integer $n\geq 0$, $\rk_{\Z}\, E(K_n) = p^n$,
$\sha(E/K_n)[p^\infty] = 0$ and $E(K_n) \otimes \zp$ is generated over
$\zp[\Gal(K_n/K)]$ by the traces of Heegner points of $p$-power conductor.

\sn
{\it Proof.\/} Theorem 4.8 and Theorem 4.9 follow from Theorem 3.4 and
Theorem 3.5, respectively, applied to $B = E$, $M = \Q$ and $\wp = p$.
Indeed, the conditions $(A1)_{E,K,p}$ and $(A5)_{E,K}$ are immediate,
$(A2)_{E,K,p}$, $(A3)_{E,K,p}$ and $(A4)_{E,K,p}$ follow from (b), (a)
and (d), respectively. Finally, $(A6)_{E,\ki/K,p}$ (resp. $(A7)_{E,\ki/K,p}$)
is a consequence of (c) and Proposition 4.6(1) (resp. of (b'), (c'), (d)
and Proposition 4.6(3)).

\sn
{\bf (4.10)} \ If $K = \Q(\sqrt{-3})$ and $p = 3$, then the
conditions (a) and (c') in Theorem 4.9 can never be satisfied simultaneously.
This is a special case of the following divisibility result, which is probably
well known, but for which we have not found any reference.

\proclaim{(4.11) Proposition}. If a prime number $p$ divides $u_K$ (i.e.,
if $(K, p) = (\Q(i), 2)$ or $(\Q(\sqrt{-3}), 3)$, then $E(K)[p] \not= 0$
or $y_K \in p E(K)$. In particural, if $y_K \not\in E(K)_{\tors}$, then
the index $[E(K) : \Z y_K]$ is divisible by $p$.

\sn
{\it Proof.\/} Assume that $E(K)[p] = 0$. According to Proposition 5.25,
there are infinitely many prime numbers $q \nmid pN$ such that
$p \nmid \widetilde E_q({\F}_q) = q + 1 - a_q$.
Any such $q$ satisfies $q \equiv \eta_K(q) \; ({\rm mod}\ {2 u_K})$,
and therefore $p \nmid (\eta_K(q) + 1 - a_q)$, since $p \mid u_K$.
The last of the norm relations in 4.4

$$
(a_q - 1 - \eta_K(q)) y_K = u_K\, \Tr_{H_q/K}(y_q) \in
u_K\, E(K) \subset p E(K)
$$
then implies that $y_K \otimes 1 \in p(E(K) \otimes \Z_{(p)})$, hence
$y_K \in p E(K)$.

\sn
{\bf (4.12)} \ It may be worthwhile to reformulate the phenomenon encountered
in Proposition~4.11 in more abstract terms, in the general situation of 4.3.
Define the {\sl group of Heegner points}

$$
E(K)_{HP} \subset E(K)
$$
to be the subgroup of $E(K)$ generated by the points

$$
y_{K,m} := \Tr_{H_m/K}(y_m),
\eqno (4.12.1)
$$
for all integers $m\geq 1$ relatively prime to $N$.

The norm relations in 4.4 imply, firstly, that $E(K)_{HP}$ is generated by
$y_{K,1} = y_K$ and by the points $y_{K,q}$ (where $q$ runs through
all primes not dividing $N$), and, secondly, that $u_K y_{K,q} \in
\Z y_K$ for all such $q$. It follows that

$$
u_K E(K)_{HP} \subset \Z y_K\subset E(K)_{HP};
$$
in particular,

$$
E(K)_{HP} = \Z y_K \qquad {\rm if\ } u_K = 1.
$$
Let us now consider the more interesting case $u_K \not= 1$, when
$(K, u_K) = (\Q(i), 2)$ or $(\Q(\sqrt{-3}, 3)$. In either case $u_K = p$
is a prime dividing $D_K$, which implies that $p\nmid N$, and therefore
$E$ has good reduction at $p$. In addition, $\chi_{p,K} = 1$ if $p = 3$
(and $\chi_{p,\Q} = 1$ if $p = 2$), where $\chi_{p,K}$ is the cyclotomic
character defined in \S{1.1}.

\proclaim{(4.13) Proposition}. Assume that $u_K = p > 1$ and $E(K)[p] \not= 0$.
\hb
(1) $E$ has good ordinary reduction at $p$, $\o\rho_{E,p} =
\pmatrix{\chi_{p, \Q} & *\cr 0 & 1\cr}$ or
$\pmatrix{1 & *\cr 0 & \chi_{p, \Q}\cr}$ in some basis of $E[p]$, and
$a_p \equiv 1 \; ({\rm mod}\ p)$.\hb
(2) For every prime $q$ not dividing $N$ we have $a_q - 1 - \eta_K(q)
\equiv 0 \; ({\rm mod}\ p)$.\hb
(3) $\Z y_K \subset E(K)_{HP} \subset \Z y_K + E(K)[p]$.

\sn
{\it Proof.\/} (1) The assumption $E(K)[p] \not= 0$ together with
$\chi_{p,K} = 1$ imply that, in a suitable basis of $E[p]$,
$\o\rho_{E,p}\vert_{G_K} = \pmatrix{1 & *\cr 0 & 1\cr}$. Therefore
$\o\rho_{E,p} = \pmatrix{\alpha\chi_{p,\Q} & *\cr 0 & \alpha\cr}$
for some character $\alpha : G_\Q \lo \Gal(K/\Q) \lo \{\pm 1\}$,
which rules out the case of supersingular reduction at $p$, by [S1, Prop. 12].

If $p = 2$, then $\alpha = 1$ and $2\nmid a_2$, for trivial reasons.
If $p = 3$, then $\alpha = 1$ or $\alpha = \chi_{3,\Q}$. In either
case, the semisimplification $\o\rho_{E,3}^{ss}$ is isomorphic to
$1 \oplus \chi_{3,\Q}$. On the other hand,
$(\o\rho_{E,3}\vert_{G_{\Q_3}})^{ss} \simeq \beta \oplus \beta \chi_{3,\Q_3}$
for an unramified character $G_{\Q_3}/I_3 \lo \{\pm 1\}$ such that
$a_3 \equiv \beta(Fr(3)) \; ({\rm mod}\ 3)$; but $\beta = 1$ by the
previous discussion.

\noindent
(2) The case $q = p$ is treated in (1). If $q \nmid pN$, then
$a_q = \Tr(\rho_{E,p}(Fr(q)) \equiv q + 1 \; ({\rm mod}\ p)$,
by (1). However, $q \equiv \eta_K(q) \; ({\rm mod}\ p)$.

\noindent
(3) For each prime $q \nmid N$, the point
$y_{K,q} - ((a_q - 1 - \eta_K(q))/p) y_K$ lies in $E(K)[p]$,
thanks to (2) and the norm relations in 4.4. In particular,
$y_{K,q} \in \Z y_K + E(K)[p]$.

\proclaim{(4.14) Proposition}. Assume that $u_K = p > 1$ and $E(K)[p] = 0$.\hb
(1) There are infinitely many primes $q \nmid pN$ satisfying
$a_q - 1 - \eta_K(q) \not \equiv 0 \; ({\rm mod}\ p)$.\hb
(2) If $q$ is as in (1), then $y_K \in p(\Z y_K + \Z y_{K,q})$ and
$E(K)_{HP} = \Z y_K + \Z y_{K,q} = \Z z_K$, where $z_K \in E(K)_{HP}$
does not depend on $q$ and satisfies $p z_K = y_K$.\hb
(3) If $y_K \in E(K)_{\tors}$ is of order $m$, then $p \nmid m$ and
$E(K)_{HP} = \Z y_K \simeq \Z/m\Z$.\hb
(4) If $y_K \not\in E(K)_{\tors}$, then $E(K)_{HP} \simeq \Z$ and
$\Z y_K = p E(K)_{HP}$.

\sn
{\it Proof.\/} (1) If $a_q - 1 - \eta_K(q) \equiv 0 \; ({\rm mod}\ p)$
for all but finitely many primes $q \nmid pN$, then $\# \widetilde E_q(\F_q)
= q + 1 - a_q \equiv q - \eta_K(q) \equiv 0 \; ({\rm mod}\ p)$ for all such $q$,
hence $E(\Q(\mu_p))[p] \not= 0$, by Proposition 5.25. This contradicts our
assumption $E(K)[p] = 0$, since $K \supset \Q(\mu_p)$.

\noindent
(2) The norm relation $p y_{K,q} = (a_q - 1 - \eta_K(q)) y_K$ together with
$p \nmid (a_q - 1 - \eta_K(q))$ imply that $y_K \in p(\Z y_K + \Z y_{K,q})$.
Fix a prime $q^\prime \nmid qN$; then $p y_{K,q^\prime} = (a_{q^\prime}
- 1 - \eta_K(q^\prime)) y_K$. There exists $n\in\Z$ such that $(a_q - 1 -
\eta_K(q)) n \equiv a_{q^\prime} - 1 - \eta_K(q^\prime) \; ({\rm mod}\ p)$;
then $p(y_{K,q^\prime} - n y_{K,q}) \in \Z p y_K$, hence
$y_{K,q^\prime} - n y_{K,q} \in \Z y_K + E(K)[p] = \Z y_K$. Therefore
$E(K)_{HP} = \Z y_K + \Z y_{K,q}$. Finally, there is a unique $z_K \in
\Z y_K + \Z y_{K,q}$ such that $p z_K = y_K$; then $y_{K,q} =
(a_q - 1 - \eta_K(q)) z_K$, hence $\Z y_K + \Z y_{K,q} = \Z z_K$.

\noindent
(3), (4) This follows from (2) and $E(K)[p] = 0$.

\sn
{\bf (4.15)} \ Let us now specialise to the case $K = \Q(\sqrt{-3})$ and
$p = 3$. Assume, in addition, that $E(K)[3] = 0$. As we saw in the proofs
of Propositions~4.11 and 4.14, there are infinitely many primes $q \nmid 3N$
such that

$$
\#\widetilde E_q(\F_q) = q + 1 - a_q \not\equiv 0 \; ({\rm mod}\ 3),
$$
thanks to Proposition~5.25 below. The point

$$
y_{K,q} = \Tr_{H_q/K}(y_q) \in E(K)
$$
satisfies

$$
3 y_{K,q} = (a_q - 1 - \eta_K(q)) y_K
$$
with $a_q - 1 - \eta_K(q) \equiv a_q - 1 - q \not\equiv 0
\; ({\rm mod}\ 3)$, and therefore $y_K \in 3 E(K)$.

Fix such a prime $q$. The discussion in \S{4.5}-\S{4.9} needs to be
modified as follows. Assume that $3 \nmid a_3$ and let, in the notation
of (4.5.1) and (4.5.2), for every integer $n\geq 0$,

$$
z_{n,q} := \alpha_3^{-n} y_{3^{n+1}q} - \alpha_3^{-n-1} y_{3^nq}
\in E(H_{3^{n+1}q}) \otimes \Z_3.
$$
These elements are again norm compatible

$$
\forall n\geq 1 \quad \Tr_{H_{3^{n+1}q}/H_{3^nq}} (z_{n,q}) =
z_{n-1,q}
$$
and the bottom element $z_{0,q} = y_{3q} - \alpha_3^{-1} y_q
\in H_{3q} \otimes \Z_3$ satisfies

$$
\displaylines{
\Tr_{H_{3q}/H_q} (z_{0,q}) = (a_3 - \sigma) y_q -
3 \alpha_3^{-1} y_q \qquad\qquad (\sigma\in \Gal(H_q/K)),\cr
\Tr_{H_{3q}/K} (z_{0,q}) = (a_3 - \sigma - \beta_3)
\Tr_{H_q/K} (y_q) = (\alpha_3 - 1) y_{K,q}.\cr}
$$

\proclaim{(4.16) Proposition}. Assume that $K = \Q(\sqrt{-3})$, $p = 3$,
$E(K)[3] = 0$ and $3 \nmid a_3$. As in 4.15, fix a prime number $q \nmid 3N$
such that $3 \nmid (a_q - 1 - q)$ and define $y_{K,q} \in E(K)$ by
(4.12.1).\hb
(1) The element $(\alpha_3 - 1) (y_{K,q} \otimes 1) \in E(K) \otimes \Z_3$
is contained in
$$
\Imm(N_{H_{3^\infty q}/K}(E \otimes \Z_3) \lo E(K) \otimes \Z_3) \subset
\Imm(N_{\ki/K}(E \otimes \Z_3) \lo E(K) \otimes \Z_3).
$$
(2) The only prime $v_3 = (\sqrt{-3})$ of $K$ above $3$ satisfies
$$
\# \widetilde E_3(k(v_3)) = \# \widetilde E_3(\F_3) \equiv 1 - \alpha_3
\; ({\rm mod}\ 3).
$$
(3) If $a_3 \not\equiv 1 \; ({\rm mod}\ 3)$, then
$3 \nmid \# \widetilde E_3(k(v_3))$ and
$$
(y_{K,q} \otimes 1) \in \Imm(N_{H_{3^\infty q}/K}(E \otimes \Z_3) \lo E(K)
\otimes \Z_3) \subset \Imm(N_{\ki/K}(E \otimes \Z_3) \lo E(K) \otimes \Z_3).
$$

\sn
{\it Proof.\/} The statements (1) and (2) follow, respectively, from the norm
relations in 4.15 and from the fact that $\# \widetilde E_3(\F_3) = 3 + 1 -
a_3$. The statement (3) is a consequence of (1) and (2).

\proclaim{(4.17) Theorem}. If $K = \Q(\sqrt{-3})$, $p = 3$ and if\hb
(a) $E(K)[3] = 0$,\hb
(b') $3\nmid a_3 \cdot (a_3 - 1) \cdot c_{\rm Tam}(E/\Q)$,\hb
(c') $y_{K,q} \not\in 3 E(K)$ (for a fixed prime $q\nmid 3N$ satisfying
$3 \nmid (a_q - 1 - q)$),\hb
(d) $\rk_{\Z}\, E(K) = 1$ and $\sha(E/K)[3^\infty] = 0$,\hb
then, for every intermediate field $K \subset L \subset \ki$,
$\sha(E/L)[3^\infty] = 0$ and the Pontryagin dual of $E(L) \otimes \Q_3/\Z_3 =
\Sel_{3^\infty}(E/L)$ is a free module of rank one over
$\Z_3[[\Gal(L/K)]]$. For every integer $n\geq 0$, $\rk_{\Z}\, E(K_n) = 3^n$,
$\sha(E/K_n)[3^\infty] = 0$ and $E(K_n) \otimes \Z_3$ is generated over
$\Z_3[\Gal(K_n/K)]$ by the traces to $K_n$ of the Heegner points of conductors
dividing $3^\infty q$.

\sn
{\it Proof.\/} The proof of Theorem 4.9 applies, except that we use
Proposition~4.16 instead of Proposition~4.6.

\vskip16pt

\centerline{\bf 5. Vanishing of certain Galois cohomology groups (after
[Ch] and [LW])}
%\bn
\vglue12pt\noindent
{\bf (5.1)} \ One of the ingredients of Kolyvagin's method for obtaining
upper bounds on the size of Selmer groups $\Sel_{p^n}(E/K) \subset
\coh 1K{E[p^n]}$ in the situation of 4.3 is a passage to the extension
$L_n := K(E[p^n])$ of $K$ over which the Galois action on $E[p^n]$ becomes
trivial. The inflation-restriction sequence

$$
0 \lo \coh 1{L_n/K}{E[p^n]} \lo \coh 1K{E[p^n]} \lo
H^1(L_n, E[p^n])^{\Gal(L_n/K)} \lo \coh 2{L_n/K}{E[p^n]}
$$
implies that such a passage entails no loss of information, provided that
$\coh 1{L_n/K}{E[p^n]} = 0$.
Sufficient criteria for the vanishing of $\coh i{L_n/K}{E[p^n]}$ were given
in [Ch, Thm. 2] (for $i = 1$); a complete answer in the case $K = \Q$ was
obtained in [LW, Thm. 1, Thm. 2] (for $i = 1, 2$). These questions were
also considered, from a slightly different point of view, in [CS1, \S{5}]
and [CS2, \S{3}].

In \S{5.2}-\S{5.21} we recall the approach adopted in [Ch] and [LW], first in
an abstract setting, then for $\wp$-power torsion in an abelian variety $B$ of
$GL(2)$-type with real multiplication (which includes the case of elliptic
curves). Unlike [Ch] and [LW], we are only interested in the ``easy case"
when $B[\wp]$ is an irreducible Galois module.

\sn
{\bf (5.2)} \ Assume that we are given the following data:

\medskip
\item{$\bullet$} \ a prime number $p$,
\item{$\bullet$} \ a finite extension $\K/\qp$, with ring of integers
$\O$, uniformiser $\pi$ and residue field $k = \O/\pi$,
\item{$\bullet$} \ a free $\O$-module $T$ of finite rank $r\geq 1$; set
$\o T := T/\pi$,
\item{$\bullet$} \ a closed subgroup $G \subset {\rm Aut}_{\O}(T) \simeq
GL_r(\O)$.

\medskip
\noindent
The $\pi$-adic filtration on $T$ induces a filtration $G = G_0 \supset G_1
\supset G_2 \supset \cdots$ by open normal subgroups

$$
G_n := \Ker(G \ho {\rm Aut}_{\O}(T) \lo {\rm Aut}_{\O}(T/\pi^n)),
$$
which have the following properties:

\medskip
\item{$\bullet$} \ $G_0/G_1 \ho {\rm Aut}_k(\o T) \simeq GL_r(k)$,
\item{$\bullet$} \ $\forall n\geq m\geq 1 \quad G_n/G_{m+n} \ho
\End_{\O}(\pi^n T/\pi^{m+n}) \simeq M_r(\O/\pi^m)$ \quad
($1 + \pi^n A \mapsto A \;\; ({\rm mod}\ \pi^m))$,
\item{$\bullet$} \ $\forall m,n\geq 1 \quad [G_m, G_n] \subset G_{m+n}$, which
implies that the adjoint action of $g \in G/G_{m+n}$ on $G_n/G_{m+n}$ (given by
${\rm ad}(g) h := g h g^{-1}$) factors through $G/G_m$.

\sn
{\bf (5.3)} \ We are interested in establishing sufficient criteria
for the vanishing of the cohomology groups $\coh 1{G/G_n}{T/\pi^m}$ (where
$n \geq m \geq 1$). Firstly, d\'evissage implies that

$$
{\rm if\ } \coh 1{G/G_n}{\o T} = 0,\ {\rm then\ }
\forall m\in \{1,\ldots,n\} \quad \coh 1{G/G_n}{T/\pi^m} = 0.
\eqno (5.3.1)
$$
Secondly, the inflation-restriction sequence for $G_n/G_{n+1} \vartriangleleft
G/G_{n+1}$ (where $n\geq 1$)

$$
0 \lo \coh 1{G/G_n}{\o T} \lo \coh 1{G/G_{n+1}}{\o T} \lo
H^1(G_n/G_{n+1}, \o T)^{G/G_n}
\eqno (5.3.2)
$$
has the following properties: $G_n/G_{n+1} \ho \End_k(\o T)$ acts
trivially on $\o T$, the action of $G/G_n$ on $\o T$ factors through
$G/G_1 \ho {\rm Aut}_k(\o T)$, and so does the adjoint action of $G/G_n$
on $\End_k(\o T)$ and its $\F_p[G/G_n]$-submodule $G_n/G_{n+1}$. As a result,

$$
H^1(G_n/G_{n+1}, \o T)^{G/G_n} = \Hom_{\F_p}(G_n/G_{n+1}, \o T)^{G/G_1}.
\eqno (5.3.3)
$$
Putting together (5.3.1)--(5.3.3), we obtain the following statement.

\proclaim{(5.4) Proposition}. Assume that $n\geq 1$ and that
$\,\forall n^\prime \in \{1,\ldots,n-1\} \quad
\Hom_{\F_p}(G_{n^\prime}/G_{n^\prime+1}, \o T)^{G/G_1} = 0$. 
If $\coh 1{G/G_1}{\o T} = 0$, then
$\,\forall m\in \{1,\ldots,n\} \quad \coh 1{G/G_n}{T/\pi^m} = 0$.

\sn
{\bf (5.5)} \ It will be convenient to investigate the conditions
in Proposition 5.4 in the following axiomatic setting. Throughout
\S{5.5}-\S{5.18},

\medskip
\item{$\bullet$} \ $p$ is a prime number,
\item{$\bullet$} \ $k$ is a finite extension of $\F_p$, of degree $f =
[k:\F_p]$,
\item{$\bullet$} \ $V$ is a finite-dimensional $k$-vector space, of
dimension $r\geq 1$,
\item{$\bullet$} \ $H \subset GL(V) \simeq GL_r(k)$ is a subgroup,
\item{$\bullet$} \ $W \subset \End_k(V) \simeq M_r(k)$ is
an $\F_p[H]$-submodule (with respect to the adjoint action of $H$).
\item{$\bullet$} \ Denote by $PH$ the image of $H$ under the projection
$GL(V) \lo PGL(V)$.

\medskip
\noindent
In order to verify the assumptions of Proposition 5.4, we must be able to
answer the following two questions (for $V = \o T$, $H = G/G_1$ and $W =
G_n/G_{n+1}$, where $n\geq 1$).

\sn
{\bf Question (Q1):} \ When is $\coh 1HV = 0$?\hb
{\bf Question (Q2):} \ When is $\Hom_{\F_p}(W, V)^H = 0$?

\medskip
There is an extensive literature devoted to (Q1); see [Gu, Thm. A]
for fairly general results (valid when $k$ is an arbitrary field
of characteristic $p$ and $H$ is a finite subgroup of $GL(V)$).

As noted in [Ch], [LW] and [CS1, CS2], one can often deduce the vanishing
statements in (Q1) and (Q2) by applying the following elementary observations.

\medskip
\item{(5.5.1)} \ If $p\nmid \# H$, then $\,\forall i > 0 \quad \coh iHV = 0$.
\item{(5.5.2)} \ (Sah's Lemma [Sa, Prop. 2.7(b)]) \ If $M$ is a $k[H]$-module
for which there
exists a central element $z\in Z(H)$ acting on $M$ by a scalar $\lambda \in
k^\times \smallsetminus \{1\}$, then $\,\forall i\geq 0 \quad \coh iHM = 0$.

\sn
{\bf (5.6)} \ Following [Ch], [LW] and [CS1, CS2], we say that $H$
{\sl contains a non-trivial homothety} if $H \cap Z(GL(V)) =
H \cap k^\times\cdot \id_V \not= \{1\}$ (or, which is equivalent,
that the projection $H \lo PH$ is not an isomorphism).

If $H$ contains a non-trivial homothety, Sah's Lemma implies that

$$
\forall i\geq 0 \quad \coh iHV = 0 = \coh iH{\Hom_{\F_p}(W, V)}.
$$
In particular, the vanishing property in both questions (Q1) and (Q2)
always holds.

\proclaim{(5.7) Proposition}. Assume that at least one of the following two
conditions is satisfied.\hb
(a) $p \nmid \# H$;\hb
(b) $V = \bigoplus V_i$ is a direct sum of simple $k[H]$-modules of
dimensions $\dim_k(V_i) \leq (p+1)/2$.\hb
Then:\hb
(1) $\End_k(V)$ is a semisimple $k[H]$-module.\hb
(2) $\End_k(V)$ is a semisimple $\F_p[H]$-module.\hb
(3) Every $\F_p[H]$-submodule $W \subset \End_k(V)$ is a direct summand.\hb
(4) If $\,\Hom_{\F_p}(\End_k(V), V)^H = 0$, then $\,\Hom_{\F_p}(W, V)^H = 0$,
for every $\F_p[H]$-submodule $W \subset \End_k(V)$.

\sn
{\it Proof.\/} The implications (2) $\Lo$ (3) $\Lo$ (4) and (a) $\Lo$ (1),(2)
are automatic, and (2) follows from (1), since the Jacobson radical of
$\F_p[H]$ is contained in the Jacobson radical of $k[H]$ ([FD, ch. 2, ex. 6,
50, 53(c)]). If $V = \bigoplus V_i$ is as in (b), so is its dual $V^* =
\bigoplus V_i^*$. Semisimplicity of the $k[H]$-module $\End_k(V) =
\bigoplus_{i, j} V_i^* \otimes V_j$ then follows from [S2, Cor. 1].

\sn
{\bf (5.8)} \ In view of Proposition 5.7, it is natural to investigate (Q2)
for $W = \End_k(V)$. In this case there is a non-degenerate $\F_p$-bilinear
symmetric pairing

$$
(\; ,\;) : W \times W \lo \F_p,\qquad (A, B) := \Tr_{k/\F_p}(\Tr(AB)),
$$
which is invariant under the adjoint action of $GL(V)$ and satisfies
$(\lambda A, B) = (A, \lambda B)$, for all $\lambda\in k$. It induces,
therefore, an isomorphism of $(k \otimes_{\F_p} k)[H]$-modules

$$
W \otimes_{\F_p} V \iso \Hom_{\F_p}(W, V).
$$
One can rewrite the tensor product on the left hand side in terms
of the Galois group

$$
\Delta := \Gal(k/\F_p) = \{\varphi^i \mid i \in \Z/f\Z\},\qquad
\varphi(a) = a^p,
$$
as follows. The ring isomorphism

$$
k \otimes_{\F_p} k \iso \prod_{\sigma\in\Delta} k,\qquad
a \otimes b \mapsto (\sigma \mapsto a \sigma(b))
$$
induces an isomorphism of $(k \otimes_{\F_p} k)[H]$-modules

$$
W \otimes_{\F_p} V \iso \bigoplus_{\sigma\in\Delta} W \otimes_k V^{(\sigma)},
\qquad V^{(\sigma)} := V \otimes_{k,\sigma} k.
$$
In concrete terms, if we fix a basis of $V$ over $k$, the (faithful) action
$\rho : H \ho GL(V) \simeq GL_r(k)$ of $H$ on $V$ gives rise to a twisted
action $\rho^{(\sigma)} : H \ho GL(V^{(\sigma)}) \simeq GL_r(k)$ given by
$\rho^{(\sigma)} = \sigma \circ \rho$.

Using this language, the $k[H]$-module $W = \End_k(V)$ corresponds
to the adjoint action $\ad(\rho) = \Hom_k(\rho, \rho) = \rho^*\otimes_k \rho
: H \lo GL(W) \simeq GL_{r^2}(k)$, and

$$
\Hom_{\F_p}(\End_k(V), V) \iso \bigoplus_{\sigma\in\Delta}
\left(\ad(\rho) \otimes_k \rho^{(\sigma)}\right).
$$
If $p \nmid \dim_k(V)$, then there is a decomposition $\ad(\rho) =
\ad^\circ(\rho) \oplus k$, where $\ad^\circ(\rho) =
\End_k^\circ(V) := \End_k(V)^{\Tr = 0}$ and the trivial representation
corresponds to the scalar endomorphisms $k \cdot \id_V$. Therefore

$$
\Hom_{\F_p}(\End_k(V), V) \iso \left(\bigoplus_{\sigma\in\Delta}
\rho^{(\sigma)}\right) \oplus \bigoplus_{\sigma\in\Delta}
\left(\ad^\circ(\rho) \otimes_k \rho^{(\sigma)}\right)
$$
if $p \nmid \dim_k(V)$. The previous discussion can be summed up as follows.

\proclaim{(5.9) Proposition}. If $\rho : H \ho GL(V)$ denotes the (faithful)
action of $H$ on $V$, then the condition $\Hom_{\F_p}(\End_k(V), V)^H = 0$
is equivalent to $\, \forall \sigma\in\Delta\quad (\ad(\rho) \otimes_k
\rho^{(\sigma)})^H = 0$. If $p \nmid \dim_k(V)$, the latter condition is
equivalent to the conjunction of $\rho^H = 0$ and $\, \forall
\sigma\in\Delta\quad (\ad^\circ(\rho) \otimes_k \rho^{(\sigma)})^H = 0$.

\sn
{\bf (5.10) A split dihedral example.} \ Assume that $p \not= 2$ and that
$n > 1$ is an odd integer dividing $\# k^\times = p^f - 1$. Denote by
$D_{2n}$ the dihedral group of order $2n$ and by $C_n \vartriangleleft D_{2n}$
its unique cyclic subgroup of order $n$. Fix an element $s\in D_{2n}
\smallsetminus C_n$; then $s^2 = 1$ and $s g s^{-1} = g^{-1}$, for all
$g\in C_n$.

For any character $\psi : C_n \lo k^\times$, the induced representation

$$
I(\psi) := {\rm Ind}_{C_n}^{D_{2n}}(\psi) : D_{2n} \lo GL(V) \simeq GL_2(k)
$$
has the following properties.

\medskip
\item{$\bullet$} \ In a suitable basis, $I(\psi)\vert_{C_n} =
\pmatrix{\psi & 0\cr 0 & \psi^{-1}\cr}$, $\;I(\psi)(s) = \pmatrix{0 & 1\cr
1 & 0\cr}$.
\item{$\bullet$} \ The image of $I(\psi)$ is contained in the normaliser
$N(C)$ of a split Cartan subgroup $C \subset GL(V)$.
\item{$\bullet$} \ $\det(I(\psi)) = \{\pm 1\} \subset k^\times$,
$\det(I(\psi)\vert_{C_n}) = \{1\}$.
\item{$\bullet$} \ $I(\psi) \simeq I(\psi)\otimes \sgn \simeq I(\psi^{-1})
\simeq I(\psi)^*$, where $\sgn : D_{2n} \lo D_{2n}/C_n \iso \{\pm 1\}$.
\item{$\bullet$} \ $I(\psi)$ is irreducible if and only if $\psi \not= 1$.
\item{$\bullet$} \ $I(1) \simeq 1 \oplus \sgn$.
\item{$\bullet$} \ $I(\psi_1) \otimes I(\psi_2) \simeq I(\psi_1\psi_2)\oplus
I(\psi_1\psi_2^{-1})$.
\item{$\bullet$} \ $\ad(I(\psi)) = I(\psi)^* \otimes I(\psi) \simeq I(\psi)
\otimes I(\psi)$, $\,\ad^\circ(I(\psi)) \simeq I(\psi^2) \oplus \sgn$.
\item{$\bullet$} \ $\forall i\in \Z/f\Z \quad I(\psi)^{(\varphi^i)} =
I(\psi^{p^i})$.
\item{$\bullet$} \ $\ad^\circ(I(\psi)) \otimes I(\psi)^{(\varphi^i)} \simeq
I(\psi^{p^i}) \oplus I(\psi^{p^i+2}) \oplus I(\psi^{p^i-2})$.
\item{$\bullet$} \ $\dim_k I(\psi)^{C_n} = 2\, \dim_k I(\psi)^{D_{2n}}$ is
equal to $2$ (resp. to $0$) if $\psi = 1$ (resp. if $\psi \not= 1$).
\medskip

\sn
{\bf (5.11) A nonsplit dihedral example.} \ Let $k_2 \simeq \F_{p^{2f}}$ be
a quadratic extension of $k$. Assume that $p \not= 2$ and that $n > 1$ is
an odd integer dividing $\# (k_2^\times/k^\times) = p^f + 1$.

For any character $\psi^\prime : C_n \lo \Ker(N_{k_2/k} : k_2^\times \lo
k^\times) \subset k_2^\times$ we define

$$
J(\psi^\prime) : D_{2n} \lo GL(V) \simeq GL_2(k)
$$
as follows. Let $V = k_2$; the regular representation $j : k_2 = \End_{k_2}(V)
\subset \End_k(V)$ identifies $k_2^\times$ with a nonsplit Cartan subgroup
$C = j(k_2^\times) \subset GL(V)$ and $\Ker(N_{k_2/k} : k_2^\times \lo
k^\times)$ with $C \cap SL(V)$. We let

$$
J(\psi^\prime)\vert_{C_n} := j \circ \psi^\prime,\qquad
J(\psi^\prime)(s) = s^\prime,
$$
for any element $s^\prime \in N(C)$ of the normaliser of $C$ with eigenvalues
$\pm 1$. Explicitly, fix $\alpha\in k_2^\times$ such that $d := \alpha^2 \in
k^\times$ and write $j$ in terms of the basis $1, \alpha$ of $k_2$ over $k$:

$$
j(a + b\alpha) = \pmatrix{a & bd\cr b & a\cr},\qquad (a, b \in k).
$$
We can then take $s^\prime = \pmatrix{1 & 0\cr 0 & -1\cr}$. The representation
$J(\psi^\prime)$ has the following properties.

\medskip
\item{$\bullet$} \ $J(\tau\circ\psi^\prime) \simeq J(\psi^\prime)$, for
any $\tau\in \Gal(k_2/k)$.
\item{$\bullet$} \ Up to isomorphism, $J(\psi^\prime)$ does not depend
on any choices.
\item{$\bullet$} \ The image of $J(\psi^\prime)$ is contained in the normaliser
$N(C)$ of a nonsplit Cartan subgroup $C \subset GL(V)$.
\item{$\bullet$} \ $\det(J(\psi^\prime)) = \{\pm 1\} \subset k^\times$,
$\det(J(\psi^\prime)\vert_{C_n}) = \{1\}$.
\item{$\bullet$} \ $J(\psi^\prime) \otimes_k k_2 \simeq I(\psi^\prime)$
(where we consider $\psi^\prime$ on the right hand side as a character
$\psi^\prime : C_n \lo k_2^\times$, and $I(\psi^\prime) : D_{2n} \lo
GL_2(k_2)$).
\medskip

\proclaim{(5.12) Proposition}. Assume that $p \not= 2$ and that $n > 1$
is an odd integer.\hb
(1) If $n \mid (p^f - 1)$ and if the character $\psi : C_n \lo k^\times$
in 5.10 is injective, then the following properties of the representation
$\rho := I(\psi) : D_{2n} \ho GL(V) \simeq GL_2(k)$ are equivalent.
$$
\displaylines{
\Hom_{\F_p}(\End_k(V), V)^{C_n} \not= 0 \iff
\Hom_{\F_p}(\End_k(V), V)^{D_{2n}} \not= 0 \iff\cr
\exists \varepsilon\in \{\pm 1\}\; \exists i \in \Z/f\Z \;\; p^i \equiv
2\varepsilon \; ({\rm mod}\ n) \iff
\exists \varepsilon\in \{\pm 1\}\; \exists i \in \Z/f\Z \;\; p^i \equiv
2\varepsilon \; ({\rm mod}\ n) {\rm\ and\ } n\mid ((2\varepsilon)^f - 1).\cr}
$$ 
(2) If $n \mid (p^f + 1)$ and if the character $\psi^\prime : C_n \lo
\Ker(N_{k_2/k})$ in 5.11 is injective, then the following properties of the
representation $\rho := J(\psi^\prime) : D_{2n} \ho GL(V) \simeq GL_2(k)$
are equivalent.
$$
\displaylines{
\Hom_{\F_p}(\End_k(V), V)^{C_n} \not= 0 \iff
\Hom_{\F_p}(\End_k(V), V)^{D_{2n}} \not= 0 \iff\cr
\exists i \in \Z/2f\Z \;\; p^i \equiv 2 \; ({\rm mod}\ n) \iff
\exists i \in \Z/2f\Z \;\; p^i \equiv 2 \; ({\rm mod}\ n) {\rm\ and\ }
n\mid (2^f - (-1)^i).\cr}
$$ 

\sn
{\it Proof.\/} The first two equivalences in (1) follow from Proposition 5.9
combined with the discussion in 5.10; the third one from the fact that the
congruences $p^i \equiv 2\varepsilon \; ({\rm mod}\ n)$ and $p^f \equiv 1
\; ({\rm mod}\ n)$ imply $(2\varepsilon)^f \equiv 1 \; ({\rm mod}\ n)$.
The statement (2) follows from the isomorphism $J(\psi^\prime) \otimes_k
k_2 \simeq I(\psi^\prime)$ combined with (1) for the pair $(\psi^\prime, k_2)$.

\sn
{\bf (5.13) Definition}
\ Given a finite field $k \simeq \F_{p^f}$ of characteristic
$p \not= 2$, we say that an odd integer $n > 1$ is {\sl $k$-exceptional}
if either $n \mid (p^f - 1)$ and the equivalent conditions in Proposition
5.12(1) are satisfied, or $n \mid (p^f + 1)$ and the equivalent conditions
in Proposition 5.12(2) are satisfied. Such a $k$-exceptional integer must
divide $2^f - 1$ or $2^f + 1$.

%\smallskip
\noindent
{\bf Examples.} (1) If $k = \F_p$, then $n$ is $k$-exceptional if and only if
$n = 3$ and $p \not= 3$.

\noindent
(2) If $k = \F_{p^2}$, then $n$ is $k$-exceptional if and only if
$n \in \{3, 5\}$ and $p \equiv \pm 2 \; ({\rm mod}\ n)$.

\noindent
(3) If $k = \F_{p^3}$, then $n$ is $k$-exceptional if and only if
$n \in \{3, 7, 9\}$ and $p \equiv \pm 2, \pm 4 \; ({\rm mod}\ n)$.

\sn
{\bf (5.14)} \ From now on, we focus our attention on the case $\dim_k(V) = 2$.
Recall Dickson's classification of subgroups $H \subset GL(V) \simeq GL_2(k)$
[Di, \S{260}].

\medskip
\item{$\bullet$} \ If $p\mid \# H$, then either $H$ acts reducibly on $V$, or
$H$ contains $SL(V^\prime)$, for some $\F_p$-vector subspace $V^\prime \subset
V$ such that $V^\prime \otimes_{\F_p} k = V$.
\item{$\bullet$} \ If $p\nmid \# H$, then either $H$ is contained in the
normaliser $N(C)$ of a Cartan subgroup $C \subset GL(V)$ (which implies that
$PH \subset PGL(V)$ is cyclic or dihedral), or $PH$ is isomorphic to
$A_4, S_4$ or $A_5$.

\medskip
\noindent
The following Proposition gives a complete list of subgroups $H \subset
GL_2(k)$ (for $p \not= 2$) acting irreducibly on $k^2$ and not containing
a non-trivial homothety (cf. [Ch, Thm. 8], [LW, Lemma 4]).

\proclaim{(5.15) Proposition}. Assume that $\dim_k(V) = 2 \not= p$. If $H
\subset GL(V)$ acts irreducibly on $V$ and does not contain a non-trivial
homothety, then:\hb
(1) There exists a Cartan subgroup $C \subset GL(V)$ such that $H \subset
N(C)$; in particular, $p\nmid \# H$.\hb
(2) The subgroup $H \cap C$ is contained in $C \cap SL(V)$; it is cyclic
of order $n > 2$, where $2 \nmid n$ and $n$ divides $\# C/\# k^\times =
\# k \mp 1$.\hb
(3) If $H \not\subset C$ (which is automatic if $C$ is split), then $H$ is
isomorphic to the dihedral group $D_{2n}$ of order $2n$, and $\det(H) =
\{\pm 1\} \subset k^\times$.\hb
[{\it In concrete terms, $H$ is isomorphic either to $D_{2n}$ or to $C_n$,
and its action on $V$ is given by $I(\psi)$ (if $H\simeq D_{2n}$) or
$J(\psi^\prime)$ (if $H\simeq D_{2n}$ or $C_n$), for an injective
character $\psi$ resp. $\psi^\prime$.\/}]\hb
(4) Conversely, if $H \subset GL(V)$ is a subgroup satisfying (1)--(2) for
some Cartan subgroup $C \subset GL(V)$ and if $H \not\subset C$ if $C$ is
split, then $H$ acts irreducibly on $V$ and does not contain a non-trivial
homothety.\hb
(5) If $k = \F_3$, then no such $H$ exists.

\sn
{\it Proof.\/} The irreducibility assumption together with the absence
of non-trivial homothety in $H$ imply, by Dickson's classification, that
$p \nmid \# H$ and that $H \simeq PH$ is cyclic, dihedral or isomorphic
to $A_4, S_4$ or $A_5$. However, the representation theory of $H$ over
$\o \F_p$ is the same as over $\C$, since $p \nmid \# H$. The groups
$A_4, S_4, A_5$ do not admit a faithful representation into $GL_2(\C)$,
therefore there is no such a representation into $GL_2(\o\F_p)$, which
leaves us only with the cases $H \simeq PH \simeq C_n$ or $D_{2n}$, for
some integer $n\geq 1$.
In particular, $H \subset N(C)$ for some Cartan subgroup $C \subset GL(V)$
and $H \cap k^\times\cdot \id_V = \{\id_V\}$, which implies that $H \cap C
\simeq P(H \cap C) \subset C/k^\times\cdot\id_V$ is cyclic of order $n > 2$
(by irreducibility), where $n \mid \# (C/k^\times\cdot\id_V)$.

If $C \simeq k_2^\times$ is nonsplit, then $n \mid (p^f + 1)$ and,
for each $a\in H\cap C$, $\det(a) = N_{k_2/k}(a) = a^{p^f + 1} = 1$.

If $C$ is split, then $n \mid (p^f - 1)$ and $H \not\subset C$. For fixed
$s\in H \smallsetminus (H \cap C)$ and any $a \in H\cap C$, $(as)^2 =
a (s a s^{-1}) = \det(a)\,\id_V \in H \cap C \cap k^\times\cdot \id_V =
\{\id_V\}$, hence $\det(a) = 1$.

In either case, the cyclic group $H \cap C$ is contained in $C \cap SL(V)$.
Its order $n > 1$ is odd, since the only element of order two in
$C \cap SL(V)$ is $- \id_V \not\in H$.

The above discussion implies that the pair $(H, \rho : H \ho GL(V))$ is of
the form $(C_n, J(\psi^\prime)\vert_{C_n}), (D_{2n}, I(\psi))$ or
$(D_{2n}, J(\psi^\prime))$, where $I(\psi)$ (resp. $J(\psi^\prime)$) is
as in 5.10 (resp. 5.11), with $\psi$ (resp. $\psi^\prime$) injective.
In each of these three cases $H$ acts irreducibly on $V$ and does not contain
a non-trivial homothety. This proves parts (1)--(4) of the proposition.
Finally, (5) follows from the fact that there is no odd $n > 2$ dividing
$3 \pm 1$.

\proclaim{(5.16) Theorem ([CS1, Thm. 9] if $k = \F_p$)}. Assume that
$\dim_k(V) = 2 \not= p$ and that $H$ acts semisimply on $V$.\hb
(1) $\forall i > 0 \quad \coh iHV = 0$.\hb
(2) If $H$ acts irreducibly on $V$, then the following conditions are
equivalent.\hb
(a) For every $\F_p$-submodule $W \subset \End_k(V)$, $\, \Hom_{\F_p}(W, V)^H
= 0$.\hb
(b) $\Hom_{\F_p}(\End_k^\circ(V), V)^H = 0$.\hb
(c) The pair $(H, \rho : H \ho GL(V))$ is not of the form
$$
(C_n, J(\psi^\prime)\vert_{C_n}),\qquad (D_{2n}, I(\psi)), \qquad
(D_{2n}, J(\psi^\prime)),
$$
for any $k$-exceptional $n > 1$.

\sn
{\it Proof.\/} (1) It is enough to assume that $p\mid \# H$, which rules out
the reducible semisimple case, when $H$ is contained in a split Cartan
subgroup. By Dickson's classification, $H$ contains $SL(V^\prime)$, which
in turn contains the homothety $-1 \in k^\times \smallsetminus \{1\}$;
we conclude by Sah's Lemma.

\noindent
(2) If $H$ contains a non-trivial homothety, (a), (b) and (c) are satisfied.
If $H$ does not contain a non-trivial homothety, then $p \nmid \# H$, by
Proposition 5.15(1). The irreducibility assumption implies that
$(V^{(\sigma)})^H = 0$, for all
$\sigma\in \Gal(k/\F_p)$. Therefore (b) is equivalent to the same statement
with $\End_k^\circ(V)$ replaced by $\End_k(V)$. The equivalence (a) $\iff$ (b)
then follows from the case (a) of Proposition 5.7, and the equivalence
(b) $\iff$ (c) from Proposition 5.12 combined with Proposition 5.15.

\proclaim{(5.17) Theorem}. In the situation of 5.2, assume that $p \not= 2
= r$, that the group $G_0/G_1$ acts irreducibly on the $k$-vector space $\o T$,
and that $G_0/G_1$ is not isomorphic to $C_n$ or $D_{2n}$, for any
$k$-exceptional odd integer $n > 1$. Then
$$
\forall m_1 \geq m_2 \geq 1 \quad \coh 1{G/G_{m_1}}{T/\pi^{m_2}T} = 0.
$$

\sn
{\it Proof.\/} Combine Theorem 5.16 with Proposition 5.4.

\proclaim{(5.18) Corollary}. Assume that $p \not= 2 = r$ and that $G_0/G_1$
acts irreducibly on the $k$-vector space $\o T$. If at least one of the
conditions (a)--(g) below holds, then
$$
\forall m_1 \geq m_2 \geq 1 \quad \coh 1{G/G_{m_1}}{T/\pi^{m_2}T} = 0.
$$
(a) $\det(G_0/G_1) \not\subset \{\pm 1\} \subset k^\times$.\hb
(b) $\det(G_0/G_1) = \{\pm 1\} \subset k^\times$ and $G_0/G_1$ is not
isomorphic to $D_{2n}$, for any $k$-exceptional odd integer $n > 1$.\hb
(c) $\det(G_0/G_1) = \{1\} \subset k^\times$ and $G_0/G_1$ is not
isomorphic to $C_n$, for any $k$-exceptional odd integer $n > 1$.\hb
(d) $\det(G_0/G_1) = {\bf F}_p^\times$ and $p > 3$.\hb
(e) $k = \F_p$ and $p = 3$.\hb
(f) $k = \F_p$, $p > 3$ and $G_0/G_1 \not\simeq A_3, S_3$.\hb
(g) $k = \F_p$ and $\det(G_0/G_1) = {\bf F}_p^\times$.

\sn
{\bf (5.19)} \ Consider the following geometric situation:

\medskip
\item{$\bullet$} \ $p \not= 2$ is a prime number,
\item{$\bullet$} \ $K$ is a field of characteristic different from $p$,
\item{$\bullet$} \ $M$ is a totally real number field,
\item{$\bullet$} \ $\wp \mid p$ is a prime of $M$ above $p$; let
$\K := M_{\wp}$, $\O := O_{\K}$, $k := \O/\wp$;
\item{$\bullet$} \ $B$ is an abelian variety over $K$ of dimension
$\dim(B) = [M : \Q]$, equipped with a ring morphism $i : O_M \lo
\End(B)$ and a symmetric $O_M$-linear isogeny $\lambda = \lambda^t :
B \lo \Bd$; let $T:= T_{\wp}(B)$ be as in (1.2).
\medskip

In this case $T$ is a free $\O$-module of rank $r = 2$. Denote by
$G \subset {\rm Aut}_{\O}(T) \simeq GL_2(\O)$ the image of the
Galois representation $\rho_{B,\wp} : G_K \lo {\rm Aut}_{\O}(T)$.
In the notation of 5.2, we have $T/\wp^n = B[\wp^n]$ (in particular,
$\o T = B[\wp]$), $G/G_n = \Gal(K(B[\wp^n])/K) \subset
{\rm Aut}_{\O}(T/\wp^n) \simeq GL_2(\O/\wp^n)$ and $G_0/G_1$ is the
image of the residual Galois representation $\o\rho_{B,\wp} : G_K \lo
{\rm Aut}_k(B[\wp]) \simeq GL_2(k)$. The Weil pairing attached to $\lambda$
implies that $\det(\rho_{B,\wp}) : G_K \lo \O^\times$ is given by the
$p$-adic cyclotomic character, hence $\det(\o\rho_{B,\wp}) = \chi_{p,K}
: G_K \lo {\bf F}_p^\times \subset k^\times$ is the $({\rm mod}\ p)$
cyclotomic character. Applying Theorem 5.17 and Corollary 5.18 in this
situation, we obtain the following results.

\proclaim{(5.20) Theorem}. In the situation of 5.19, assume that
$B[\wp]$ is an irreducible $k[G_K]$-module, and that $\o\rho_{B,\wp}(G_K)
\subset GL_2(k)$ is not isomorphic to $C_n$ or $D_{2n}$, for any
$k$-exceptional odd integer $n > 1$. Then
$$
\forall m_1 \geq m_2 \geq 1 \quad \coh 1{K(B[\wp^{m_1}])/K}{B[\wp^{m_2}]} = 0.
$$

\proclaim{(5.21) Corollary}. Assume that $B[\wp]$ is an irreducible
$k[G_K]$-module. If at least one of the conditions (a)--(g') below holds, then
$$
\forall m_1 \geq m_2 \geq 1 \quad \coh 1{K(B[\wp^{m_1}])/K}{B[\wp^{m_2}]} = 0.
$$
(a) $\chi_{p,K}(G_K) \not\subset \{\pm 1\} \subset {\bf F}_p^\times$.\hb
(a') $K \supset \Q$ and $\Q(\mu_p)^+ \not \subset K$.\hb
(a'') $K$ is an imaginary quadratic field and $p > 3$.\hb
(b) $\chi_{p,K}(G_K) = \{\pm 1\} \subset {\bf F}_p^\times$ and
$\o\rho_{B,\wp}(G_K)$ is not isomorphic to $D_{2n}$, for any $k$-exceptional
odd integer $n > 1$.\hb
(b') $\Q(\mu_p)^+ \subset K$, $\Q(\mu_p)\not\subset K$ and
$\o\rho_{B,\wp}(G_K)$ is not isomorphic to $D_{2n}$, for any $k$-exceptional
odd integer $n > 1$.\hb
(c) $\chi_{p,K}(G_K) = \{1\} \subset {\bf F}_p^\times$ and
$\o\rho_{B,\wp}(G_K)$ is not isomorphic to $C_n$, for any $k$-exceptional
odd integer $n > 1$.\hb
(c') $\Q(\mu_p) \subset K$ and $\o\rho_{B,\wp}(G_K)$ is not isomorphic to
$C_n$, for any $k$-exceptional odd integer $n > 1$.\hb
(d) $\chi_{p,K}(G_K) = {\bf F}_p^\times$ and $p > 3$.\hb
(d') $K \supset \Q$, $K \cap \Q(\mu_p) = \Q$ and $p > 3$.\hb
(e) $k = \F_p$ and $p = 3$.\hb
(e') $K$ is an imaginary quadratic field and $k = \F_p$.\hb
(f) $k = \F_p$, $p > 3$ and $\o\rho_{B,\wp}(G_K) \not\simeq A_3, S_3$.\hb
(g) $k = \F_p$ and $\chi_{p,K}(G_K) = {\bf F}_p^\times$.\hb
(g') $K \supset \Q$, $K \cap \Q(\mu_p) = \Q$ and $k = \F_p$.

\sn
{\bf (5.22)} \ If $M = \Q$, then $\wp = p$, $\K = \qp$, $\O = \zp$ and
$B = E$ is an elliptic curve. In this case much more precise results
were proved in [Ch, Thm. 2] and [LW, Thm. 11], under suitable assumptions
on $K$.

We now prove several auxiliary results that will be needed in \S{6}
(Proposition~5.25 was already used in the proofs of Propositions~4.11
and 4.14).

\proclaim{(5.23) Proposition}. Let $V$ be a two-dimensional vector space
over a field $k$.\hb
(1) If $G \subset GL(V)$ is a subgroup satisfying
$\, \forall g\in G \quad \det(1 - g \mid V) = 0$, then there exists
a one-dimensional subspace $W \subset V$ such that $W^G \not= 0$ or
$(V/W)^G \not= 0$. Equivalently, there exists a basis of $V$ in which
$G \subset H_1 := \pmatrix{1 & *\cr 0 & *\cr}$ or $G \subset
H_2 := \pmatrix{* & *\cr 0 & 1\cr}$.\hb
(2) If $G \subset GL(V)$ is a subgroup satisfying
$\, \forall g\in G \quad \Tr(g - 1 \mid V) = 0$ and if the characteristic
of $k$ is not equal to $2$, then there exists a one-dimensional subspace
$W \subset V$ such that $W^G = W$ and $(V/W)^G = V/W$. Equivalently,
there exists a basis of $V$ in which $G \subset \pmatrix{1 & *\cr 0 & 1\cr}$.

\sn
{\it Proof.\/} (1) The eigenvalues of any $g\in G$ are equal to $1$ and
$\det(g)$.  In particular, if $g \in G \cap SL(V)$, then $g$ is unipotent
and $\Tr(g) = 2$, $\det(g) = 1$.

If $\# (G \cap SL(V)) > 1$, then there exists a basis of $V$ in which $g_0
:= \pmatrix{1 & 1\cr 0 & 1\cr} \in G$. If $g = \pmatrix{a & b\cr c & d\cr}
\in G \cap SL(V)$, then $ad - bc = 1$ and $a + d = \Tr(g) = 2 = \Tr(g_0 g) =
a + c + d$, which implies that $c = 0$, and both eigenvalues of $g$ are equal
to $a = d = 1$; thus $G \cap SL(V) \subset H_1 \cap H_2$ and
$G \subset \{g \in GL(V) \mid g g_0 g^{-1} \subset H_1 \cap H_2\} =
\pmatrix{* & *\cr 0 & *\cr}$. This means that $G$ is contained in the union
of the subgroups $H_1$ and $H_2$ of $GL(V)$, and therefore is contained
in one of them.

If the group $G \cap SL(V)$ is trivial, then $\det : G \iso \det(G) \subset
k^\times$ is an isomorphism and $G$ is abelian. As a result, for each
$g\in G \smallsetminus \{\id_V\}$, the direct sum decomposition
$V = V^{g=1} \oplus V^{g = \det(g)}$ is $G$-stable, hence
$G \subset H_1^\prime \cup H_2^\prime$, where $H_1^\prime :=
\pmatrix{1 & 0\cr 0 & *\cr}$ and $H_2^\prime := \pmatrix{* & 0\cr 0 & 1\cr}$.
Again, this implies that $G \subset H_1^\prime$ or $G \subset H_2^\prime$.

\noindent
(2) For each $g\in G$ we have $2 \det(1 - g \mid V) = \Tr(g - 1) \Tr(g) -
\Tr(g^2 - 1) = 0$. Part (1) then implies that there exists $i \in \{1, 2\}$
such that $G \subset H_i^{\Tr = 2} = H_1 \cap H_2$.

\proclaim{(5.24) Corollary}. Assume that, in the situation of 5.2, $r = 2$
and $Y \subset G$ is a subset that maps surjectively on $G_0/G_1$.\hb
(1) If $\,\forall g \in Y \;\det(1 - g \mid T) \equiv 0 \; ({\rm mod}\ \pi)$,
then there is a basis of $\o T$ in which $G_0/G_1 \subset
\pmatrix{1 & *\cr 0 & *\cr}$ or $G_0/G_1 \subset \pmatrix{* & *\cr 0 & 1\cr}$.
\hb
(2) If $p\not= 2$ and $\,\forall g \in Y \;\Tr(g - 1 \mid T) \equiv 0
\; ({\rm mod}\ \pi)$, then there is a basis of $\o T$ in which
$G_0/G_1 \subset \pmatrix{1 & *\cr 0 & 1\cr}$.

\proclaim{(5.25) Proposition}. If, in the situation of 5.19, $K$ is a number
field and there exists a finite set $S$ of finite primes of $K$ (containing
all primes above $p$ and all primes at which $B$ has bad reduction) such that
$$
\forall v\not\in S \quad \# \widetilde B_v(k(v)) \equiv 0 \; ({\rm mod}\ p),
$$
then $\o\rho_{B,\wp}$ is isomorphic to $\pmatrix{1 & *\cr 0 & \chi_{p,K}\cr}$
or $\pmatrix{\chi_{p,K} & *\cr 0 & 1\cr}$. In particular, $B(K(\mu_p))[\wp]
\not= 0$.

\sn
{\it Proof.\/} For each $v\not\in S$,

$$ {\rm det}_{\O}(1 - Fr(v) \mid T_{\wp}(B)) =
\# \widetilde B_v(k(v)) \equiv 0 \; ({\rm mod}\ \wp).
$$
The statement of the proposition follows from Corollary 5.24(1) applied to
$T = T_{\wp}(B)$, $G = \Imm(G_K \lo {\rm Aut}_{\O}(T))$ and $Y =
\{Fr(v) \mid v\not\in S\}$ (which maps surjectively on $G_0/G_1 =
\Imm(G_K \lo {\rm Aut}_k(B[\wp]))$, by the \v{C}ebotarev density theorem
for $K(B[\wp])/K$).

\proclaim{(5.26) Proposition}. Let $E$ be an elliptic curve over $\Q$
of conductor $N$ and $K$ a quadratic field of discri\-mi\-nant $D_K$ relatively
prime to $N$.  Let $\rho := \o\rho_{E,p} : G_{\Q} \lo {\rm Aut}_{\F_p}(E[p])
\simeq GL_2(\F_p)$, for a prime number $p \not= 2$.\hb
(1) The field $L := \Q(E[p])$ has the following property:
$$
\rho(G_{\Q}) \not= \rho(G_K) \; \iff\; L\cap K = K \;\iff\; D_K = p^*
:= (-1)^{(p-1)/2} p.
$$
(2) If $\rho$ is irreducible, so is $\rho\vert_{G_K}$.\hb
(3) If $\rho\vert_{G_K}$ is irreducible, but not absolutely irreducible,
then $p = 3$, $K = \Q(\sqrt{-3})$, $E$ has good ordinary reduction at $3$,
$\rho(G_K)$ is a cyclic group of order $4$ and $\rho(G_{\Q})$ is a dihedral
group of order $8$.

\sn
{\it Proof.\/} (1) We have $\rho(G_{\Q}) \iso \Gal(L/\Q)$ and $\Gal(L/L\cap K)
\simeq \Gal(KL/K) \iso \rho(G_K)$, which yields the first equivalence in (1).
A prime number $\ell$ is unramified in $K/\Q$ if and only if $\ell \nmid D_K$;
it is unramified in $L/\Q$ if $\ell \nmid pN$. As $(N, D_K) = 1$, the equality
$L\cap K = K$ implies that $\{\ell \mid D_K\} \subset \{\ell \mid D_K\} \cap
\{\ell \mid pN\} \subset \{p\}$, hence $D_K = p^*$. Conversely,
$\Q(\sqrt{p^*}) \subset \Q(\mu_p) \subset L$.

\noindent
(2) If $\rho$ is irreducible but $\rho\vert_{G_K}$ is not, then
$\rho\vert_{G_K}$ is semisimple (since $G_K$ is a normal subgroup of $G_{\Q}$)
and its image is contained in a split Cartan subgroup $C_s$ of $GL_2(\F_p)$.
Moreover, $\rho(G_{\Q}) \not= \rho(G_K)$, hence $D_K = p^*$ and $p \nmid N$,
which means that $E$ has good reduction at $p$.

If the reduction at $p$ is supersingular, then $\rho(G_{\qp}) = N(C_{ns})$
is the normaliser of a nonsplit Cartan subgroup $C_{ns}$ of $GL_2(\F_p)$,
by [S1, Prop. 12]. In particular, $\# \rho(G_K)$ is a multiple of
$\# N(C_{ns})/2 = p^2 - 1 > (p - 1)^2 = \# C_s \geq \# \rho(G_K)$,
which is impossible.

If the reduction at $p$ is ordinary, then the restriction of $\rho$ to
the inertia group $I_p \subset G_{\qp}$ is given by
$\pmatrix{\chi_{p, \qp} & *\cr 0 & 1\cr}$, by [S1, Prop. 11]. On the other
hand, $\rho\vert_{G_K} = \alpha \oplus \alpha^c$, where $\alpha : G_K \lo
{\bf F}_p^\times$ is a character and $\alpha^c(g) := \alpha(\tilde c g
\tilde c^{-1})$, for any $\tilde c \in G_{\Q} \smallsetminus G_K$.
Consequently, the restrictions to the inertia group $I_{\wp} \subset
G_{K_{\wp}}$ (where $\wp \mid p$ is the only prime of $K$ above $p$)
satisfy $\{\alpha\vert_{I_{\wp}}, \alpha^c\vert_{I_{\wp}}\} =
\{\chi_{p, K_{\wp}}\vert_{I_{\wp}}, 1\}$. As a result,
$\chi_{p, K_{\wp}}\vert_{I_{\wp}} = 1$, which implies that
$\chi_{p, \qp}^2(I_p) = 1$, $p = 3$ and $K = \Q(\sqrt{-3})$. In this case
$\chi_{3, K} = 1$, hence $\rho(G_K) \subset C_s \cap SL_2(\F_3) = \{\pm I\}$.
As $\rho(G_{\Q})$ contains $\rho(\tilde c) \sim \pmatrix{1 & 0\cr 0 & -1\cr}$,
we have $\rho(G_{\Q}) \simeq (\Z/2\Z)^a$ for some $a \leq 2$, which contradicts
the irreducibility of $\rho$.

\noindent
(3) Firstly, $\rho(G_K)$ is contained in $C_{ns}$ but not in $C_{ns} \cap
C_s = {\bf F}_p^\times \cdot I$. Secondly, $\rho(G_{\Q})$ contains
$\rho(\tilde c) \not\in C_{ns}$, hence $\rho(G_{\Q}) \not= \rho(G_K)$;
thus $D_K = p^*$ and $E$ has good reduction at $p$.

If the reduction at $p$ is supersingular, then

$$
\# \rho(G_K) = \# \rho(G_{\Q})/2 \geq \# \rho(G_{\qp})/2 = \# C_{ns} \geq
\# \rho(G_K).
$$
It follows that $\rho(G_K) = C_{ns}$ and $\det \rho(G_K) =
N_{\F_{p^2}/\F_p}({\bf F}_{p^2}^\times) = {\bf F}_p^\times$, which is
equivalent to $\Q(\mu_p) \cap K = \Q$, but this is not true.

If the reduction at $p$ is ordinary, then the restriction of $\rho$ to
$I_p$ is of the form $\pmatrix{\chi_{p, \qp} & *\cr 0 & 1\cr} \subset C_{ns}$,
which implies again that $\chi_{p, K_{\wp}}\vert_{I_{\wp}} = 1$, $p = 3$,
$K = \Q(\sqrt{-3})$ and $\chi_{3, K} = 1$, hence $\rho(G_K)$ is contained in
$C_{ns} \cap SL_2(\F_3)$, which is a cyclic group of order $4$. On the
other hand, $\# \rho(G_K) > 2$, by the irreducibility of $\rho\vert_{G_K}$,
which implies the statements about the structure of $\rho(G_K)$ and
$\rho(G_{\Q})$.

\sn
{\bf (5.27) Genus theory of quadratic fields.} \ Let $K$ be a quadratic
field, $R = \{q \mid D_K\}$ the set of prime numbers ramified in $K/\Q$,
$C$ the strict ideal class group of $K$, $H$ the strict Hilbert class
field of $K$ (the maximal abelian extension of $K$ unramified over $K$
at all finite primes) and $K_{\gen} := H \cap \Q^{ab}$ the genus field
of $K$. The Galois groups in the tower $\Q \ho K \ho K_{\gen} \ho H$
are as follows.

$$
\displaylines{
G := \Gal(H/K) \simeq C,\qquad G_+ := \Gal(H/\Q) {\rm\ \ satisfies\ \ }
\forall g_+ \in G_+\smallsetminus G\;\; \forall g \in G \;\; g_+^2 \in G,
\quad g_+ g g_+^{-1} = g^{-1}\cr
\Gal(H/K_{\gen}) = [G_+, G_+] = G^2 \simeq C^2,\qquad \Gal(K_{\gen}/K)
\simeq C/C^2.\cr}
$$
There is a unique factorisation

$$
D_K = \prod_{q\in R} D_q,\qquad D_q \equiv 0, 1 \; ({\rm mod}\ 4),\qquad
|D_q| = {\rm\ a\ power\ of\ } q
$$
(if $q \not= 2$, then $D_q = q^* := (-1)^{(q-1)/2} q$). In terms of this
factorisation,

$$
K_{\gen} = \Q(\{\sqrt{D_q}\}_{q\in R})
$$
is the compositum of the quadratic fields $K(q) := \Q(\sqrt{D_q})$, for all
$q \in R$.

\proclaim{(5.28) Proposition}. For each $q\in R$, the compositum $H(q)$
of all subfields of $H$ unramified over $\Q$ outside $q \infty$ is equal to
$$
H(q) = \cases{
H, & if $R = \{q\}$\cr
K(q), & if $R \not= \{q\}$.\cr}
$$

\sn
{\it Proof.\/} The case $R = \{q\}$ is immediate. Assume that $R \not= \{q\}$.
For each $q^\prime \in R \smallsetminus \{q\}$ and each prime $v$ in $H$
above $q^\prime$, the inertia subgroup $I_v \subset \Gal(H/\Q) = G_+$ is of the
form $I_v = \{1, h_v\}$, where $h_v^2 = 1$ and $h_v \not= G$. By definition,
$H(q)$ is the fixed field of the subgroup $G(q) \subset G$ generated by the
$I_v$, for all $q^\prime \in R \smallsetminus \{q\}$ and $v \mid q^\prime$.
If $g\in G$, then $g h_v g^{-1} \in I_{g(v)}$ and $g^2 = g h_v g^{-1} h_v^{-1}
\in G(q)$; thus $G^2 \subset G(q)$ and $H(q) \subset K_{\gen}$, but the only
subfields of $K_{\gen}$ unramified over $\Q$ outside $q \infty$ are $\Q$
and $K(q)$.

\sn
{\bf (5.29)} \ For an arbitrary quadratic field $K$, its ring class field
$H_n$ of conductor $n \geq 1$ is an abelian extension of $K$ characterised
by the fact that the reciprocity map of class field theory induces
an isomorphism

$$
K^\times_+\backslash \widehat K^\times/\widehat O_n^\times \iso \Gal(H_n/K),
$$
where $\widehat K = K \otimes \widehat \Z$, $\widehat O_n = (\Z + n O_K)
\otimes \widehat \Z$ and $K^\times_+ \subset K^\times$ is the subgroup of
elements that are positive under all real embeddings $K \ho \R$. For $n = 1$,
$H_1$ is the strict Hilbert class field of $K$. In general, $H_n$ is a Galois
extension of $\Q$ and

$$
\forall g\in \Gal(H_n/K) \quad \forall g_+\in \Gal(H_n/\Q)\smallsetminus
\Gal(H_n/K) \qquad
g_+^2 \in \Gal(H_n/K), \quad g_+ g g_+^{-1} = g^{-1}.
$$
In particular, $\Gal(H_n/\Q)^{ab} \iso (\Z/2\Z)^a$ for some $a\geq 0$.

\sn
{\bf (5.30)} \ In the situation of 5.2, assume that we are given surjective
morphisms $G_{\bf \Q} \, \mapr\rho\, G \, \mapr \chi\, {\bf Z}_p^\times$
whose composition is the cyclotomic character, and a surjective $\O$-bilinear
pairing $\langle\; ,\; \rangle : T \times T \lo \O$ satisfying

$$
\forall g\in G \quad \forall x,y\in T \qquad \langle gx, gy \rangle =
\chi(g) \langle x, y \rangle.
$$
For each $m\geq 1$, let $\rho_m$ be the composition $\rho_m : G_{\Q} \lo G
\lo G/G_m \ho {\rm Aut}_{\O/\pi^m}(T/\pi^m)$ and define $L_m := \Q(T/\pi^m T)
= \qb^{\Ker(\rho_m)}$.

By definition, if $g \in G_m$ ($m \geq 1$), then $(g - 1)T \subset \pi^m T$
and

$$
\forall x,y\in T \quad (\chi(g) - 1) \langle x, y\rangle =
\langle gx, gy \rangle - \langle x, y\rangle = \langle (g - 1)x, gy\rangle
+ \langle x, (g - 1)y\rangle \in \pi^m \O,
$$
hence $\chi(g) \in 1 + \pi^m\O$, by the surjectivity of $\langle\; ,\;
\rangle$. This implies that

$$
\forall m\geq 1 \quad L_m = \Q(T/\pi^m T) \supset \Q(\mu_{p^t}),
$$
where $t$ is the smallest integer such that $t \geq m/e$ and $e :=
\ord_{\pi}(p)$ is the ramification index of $\K/\qp$. In particular,

$$
L_\infty := \bigcup_{m\geq 1} L_m \supset \Q(\mu_{p^\infty}).
$$

\proclaim{(5.31) Proposition}. Assume that we are in the situation of 5.30
with $p \not= 2$, that $K$ is a quadratic field of discriminant $D_K$ and that
$\rho : G_{\Q} \lo {\rm Aut}_{\O}(T)$ is unramified outside $pN\infty$ (i.e.,
that $L_\infty/K$ is unramified outside $pN\infty$). Fix $m, n \geq 1$.\hb
(1) For every algebraic extension $F/\Q$ we have $(T/\pi^m T)^{G_F} =
(T/\pi^m T)^{G_{F\cap L_m}}$.\hb
(2) If $L_m \subset H_n$, then $p = 3$, $1 \leq m \leq e$ and $3 \mid
n D_K$.\hb
(3) If $(n, pN) = 1$, then $K L_\infty \cap H_n = K L_\infty \cap H_1$.\hb
(4) If $(N, D_K) = 1$, then the extension $(L_\infty \cap H_1)/\Q$
is unramified outside $p \infty$.\hb
(5) If $(pN, D_K) = 1$, then $L_\infty \cap H_1 = \Q$.\hb
(6) If $(N, D_K) = 1$ and $D_K = p^* := (-1)^{(p-1)/2}$, then $K \subset
L_1$.\hb
(7) If $(N, D_K) = 1$, $p \mid D_K$ and $D_K \not= p^*$, then $L_\infty \cap
H_1 = \Q(\sqrt{p^*}) = L_1 \cap H_1$ and $L_\infty \cap K = \Q$.\hb
(8) If $(N, D_K) = 1$, $D_K = p^*$ and $r = \rk_{\O} (T) = 2$, then
$\o T^{G_{H_1}} = \o T^{G_K}$.

\sn
{\it Proof.\/} (1) This is true by the definition of $L_m$.

\noindent
(2) If $t$ is the smallest integer such that $t \geq m/e$,
then $L_m \subset H_n$ implies that $\Q(\mu_{p^t}) \subset L_m \cap \Q^{ab}
\subset H_n \cap \Q^{ab} =$ a compositum of quadratic fields unramified
outside $n D_K \infty$. Therefore $\varphi(p^t) \leq 2$, $p = 3$, $t = 1$
and $3 \mid n D_K$.

\noindent
(3) The extension $K L_\infty/K$ (resp. $H_n/K$) is unramified
outside $\{v \mid pN\infty\}$ (resp. $\{v \mid n\infty\}$);
thus $(KL_\infty \cap H_n)/K$ is an abelian extension unramified
at all finite places, so it must be contained in $H_1$.

\noindent
(4) The extension $L_\infty/\Q$ (resp. $H_1/\Q$) is unramified
outside $\{\ell \mid pN\infty\}$ (resp. $\{\ell \mid D_K\infty\}$);
thus $(L_\infty \cap H_1)/\Q$ is unramified outside $p \infty$.

\noindent
(5) In this case $(L_\infty \cap H_1)/\Q$ is unramified outside $\infty$,
so $L_\infty \cap H_1 = \Q$.

\noindent
(6) $K = K(p) := \Q(\sqrt{p^*}) \subset \Q(\mu_p) \subset L_1$.

\noindent
(7) The quadratic field $K(p) = \Q(\sqrt{p^*})$ is contained in both
$\Q(\mu_p) \subset L_1$ and in $H_1$ (by genus theory); thus $K(p) \subset
L_1 \cap H_1 \subset L_\infty \cap H_1$. On the other hand, (4) tells us that
$L_\infty \cap H_1$ is contained in $H(p)$, but $H(p) = K(p)$ in our case,
by Proposition 5.28.

\noindent
(8) If $p = 3$, then $K = \Q(\sqrt{-3}) = H_1$.
If $p > 3$, then $L_1 \not\subset H_1$ by (2), which means that
$d := \dim_k \o T^{G_{H_1}} \leq 1$. There is nothing to prove if $d = 0$.
If $d = 1$, then $\Gal(L_1 \cap H_1/\Q)$ acts on the line $\o T^{G_{H_1}} =
\o T^{G_{L_1\cap H_1}}$ by a character $\alpha : \Gal(L_1 \cap H_1/\Q) \lo
\Gal(L_1 \cap H_1/\Q)^{ab} \lo k^\times$. However, $\Gal(L_1 \cap H_1/\Q)^{ab}$
is a quotient of $\Gal(H_1/\Q)^{ab} = \Gal(K_{\rm gen}/\Q) = \Gal(K/\Q)$,
which means that $G_K$ acts on $\o T$ by $\pmatrix{1 & *\cr 0 & \chi_{p,K}\cr}$.
As $\chi_{p,K} \not= 1$ for $p > 3$, it follows that $\o T^{G_K} =
\o T^{G_{H_1}}$, as claimed.

\vskip16pt

\centerline{\bf 6. Kolyvagin's result on the vanishing of $\sha(E/K)[p^\infty]$}

\bn
{\bf (6.1)} \ Throughout \S{6}, let:

\medskip
\item{$\bullet$} \quad $E$ be an elliptic curve over $\Q$ of conductor $N$,
\item{$\bullet$} \quad $\varphi : X_0(N) \lo E$ a modular parameterisation
of $E$ sending $i\infty$ to the origin,
\item{$\bullet$} \quad $K$ an imaginary quadratic field in which all primes
dividing $N$ split,
\item{$\bullet$} \quad ${\cal N}$ an ideal of $O_K$ such that $O_K/{\cal N}
\simeq \Z/N\Z$.

\medskip
\noindent
As in 4.3, these data determine the Heegner points $y_m \in E(H_m)$ on $E$,
defined over the ring class fields $H_m$ of conductors $m\geq 1$ relatively
prime to $N$, and the basic Heegner point $y_K = \Tr_{H_1/K}(y_1)$.

\sn
{\bf (6.2)} \ If $y_K \not\in E(K)_{\tors}$ (and $D_K \not= -3, -4$), then
the groups $E(K)/\Z y_K$ and $\sha(E/K)$ are finite ([K1, Thm. A]) and the
N\'eron--Tate height of $y_K$ is given by the formula of Gross and Zagier
[GZ, Thm. V.2.1] (Gross and Zagier considered only the case when $D_K$ is odd;
for even $D_K$ the corresponding formula is a special case of [Z, Thm. 1.2.1]).
Combining their formula with the conjecture of Birch and Swinnerton-Dyer,
Gross and Zagier observed [GZ, Conj. V.2.2] that, if $y_K \not\in
E(K)_{\tors}$, then the conjecture of Birch and Swinnerton-Dyer for $E$
over $K$ holds if and only if

$$
[E(K) : \Z y_K] \buildrel ? \over = (\# \sha(E/K))^{1/2}\, u_K\,
c_{\rm Tam}(E/\Q)\, c_{\rm Manin}(\varphi),
\eqno (6.2.1)
$$
where $c_{\rm Tam}(E/\Q) = \prod_{\ell\mid N} c_{{\rm Tam},\ell}(E/\Q)$ is the
product of all non-archimedean local Tamagawa factors of $E$ over $\Q$,
$u_K = \# (O_K^\times/\Z^\times)$ and $c_{\rm Manin}(\varphi) \in \Z_{> 0}$
is the Manin constant for $\varphi$.

Recall that, for any elliptic curve $E^\prime$ defined over any number field
$K^\prime$, the Cassels--Tate pairing on the finite abelian group
$\sha(E^\prime/K^\prime)/\divv$ with values in $\Q/\Z$ is alternating and
non-degenerate, which implies that $\sha(E^\prime/K^\prime)/\divv$ is of the
form $X \oplus X$, for some maximal isotropic subspace $X$. In particular,
$\# (\sha(E^\prime/K^\prime)/\divv) = (\# X)^2$ is a square.

In \S{0.8}-\S{0.9} we discussed Kolyvagin's results on a conjectural
divisibility

$$
{\rm if\ } y_K \not\in E(K)_{\tors}, {\rm\ then} \quad
[E(K) : \Z y_K]/(\# \sha(E/K))^{1/2} \in \Z_{(p)},
\eqno (6.2.2)
$$
for a fixed prime $p \not= 2$.
Jetchev [J, Thm. 1.1] proved, under suitable assumptions, a sharpening of
(6.2.2) in the following form:

$$
{\rm if\ } y_K \not\in E(K)_{\tors}, {\rm\ then} \quad \forall \ell\mid N\quad
[E(K) : \Z y_K]/((\# \sha(E/K))^{1/2} c_{{\rm Tam},\ell}(E/\Q)) \in \Z_{(p)},
\eqno (6.2.3)
$$
in line with (6.2.1).

\sn
{\bf (6.3)} \ The simplest case of the expected divisibility (6.2.2) is the
following statement:

$$
{\rm if\ } y_K \not\in p E(K) + E(K)_{\tors}, {\rm\ then\ } E(K) \otimes\zp
= \zp(y_K \otimes 1) \simeq \zp {\rm\ and\ } \sha(E/K)[p^\infty] = 0
\eqno (6.3.1)
$$
(if $E(K)[p] = 0$, then $p E(K) + E(K)_{\tors} = p E(K)$). As recalled
in \S{0.3}-\S{0.4}, (6.3.1) was deduced by Kolyvagin [K1] from his more
general annihilation result [K1, Cor. 13] under the assumption that
$p \not= 2$, $u_K = 1$ and $\rho := \o\rho_{E,p} : G_{\Q} \lo
{\rm Aut}_{\F_p}(E[p]) \simeq GL_2(\F_p)$ has ``large image".

A more direct exposition of Kolyvagin's proof of (6.3.1) in the case when
$p \nmid 2 D_K$ and $\rho$ is surjective was given by Gross
[G, Proposition 2.1, Proposition 2.3]. It turns out that the arguments in [G]
are valid under weaker assumptions, as we are now going to explain. We begin
by extracting from [G] the conditions on $E$, $K$ and $p$ used in the proof.
After that we show that only one of them (an irreducibility assumption)
really matters.

\proclaim{(6.4) Proposition ([G, Prop. 2.1, Prop. 2.3 and its proof]}.
If $p\not= 2$ is a prime number and if the conditions (C1)--(C6) below
are satisfied, then the implication (6.3.1) holds.\hb
(C1) $u_K = 1$ (i.e., $D_K \not= -3, -4$).\hb
(C2) For each $n\geq 1$ relatively prime to $pND_K$, $E(H_n)[p] = 0$.\hb
(C3) $E(\Q)[p] = 0$.\hb
(C4) For $i = 1,2$,\quad $H^i(K(E[p])/K, E[p]) = 0$.\hb
(C5) The restriction of $\rho = \o\rho_{E,p}$ to $G_K$ is absolutely
irreducible.\hb
(C6) Neither of the two subgroups $E[p]^\pm \subset E[p]$ (:= the
$(\pm 1)$-eigenspaces for the action of complex conjugation) contains
a non-zero $G_K$-stable subgroup (equivalently, $\rho$ is irreducible).

\sn
{\it Proof.\/} The conditions (C1) and (C2) are used in [G, \S{3}-\S{5}]
in order to construct Kolyvagin's derivative classes and establish their
basic properties, and (C3) is needed in the proof of [G, Prop. 6.2(1)]
for $v \mid N$. In the general discussion in [G, \S{7}-\S{8}], no additional
conditions are needed. Things begin to get more interesting in [G, \S{9}].
The condition (C4) implies the statement of [G, Prop. 9.1] (the proof of
which relied on the assumption that $p \nmid D_K$; this was not stated
explicitly in [G, Prop. 2.1, Prop. 2.3]). The irreducibility conditions
(C5) and (C6) are used, respectively, in the proofs of [G, Prop. 9.3]
and [G, Prop. 9.5(2)] (note that the proof of [G, Prop. 9.3] relies on the
fact that $\End_{\F_p[G_K]}(E[p]) = \F_p$). The rest of the proof in
[G, \S{9}-\S{10}] goes through unchanged.

\proclaim{(6.5) Proposition}. For any prime number $p \not= 2$, the
conditions (C2), (C3), (C4) and (C6) in Proposition~6.4 follow from (C5).
Therefore the implication (6.3.1) holds if $p \not= 2$, $D_K \not= -3, -4$
and $E[p]$ is an absolutely irreducible $\F_p[G_K]$-module (the latter
condition implies that $E(K)[p] = 0$).

\sn
{\it Proof.\/} The implication (C5) $\Lo$ (C6) is straightforward, and
(C3) follows from (C6) and the fact that $\dim_{\F_p} E(\Q)[p] \leq 1$.
The implication (C5) $\Lo$ (C2) is a special case of Proposition 5.31(8).
Finally, (C4) follows from Sah's Lemma (5.5.2) and the fact that
$\rho(G_K) \subset GL_2(\F_p)$ contains a non-trivial homothety
(by Proposition 5.15, since $\# \det(\rho(G_K)) = \# \chi_{p,K}(G_K)
> 2$ for $p > 3$).

\sn
{\bf (6.6)} \ We are now ready to reprove (and slightly extend) the refinement
of Kolyvagin's result on (6.3.1) established by Cha [Ch, the case $m = 0$
of Thm. 21].

\proclaim{(6.7) Theorem}. Assume that $p \not= 2$ and that $E[p]$ is an
irreducible $\F_p[G_{\Q}]$-module (which implies that $E(K)[p] = 0$).\hb
(1) Assume that $(K, p) \not= (\Q(\sqrt{-3}), 3)$. If $y_K \not\in p E(K)$,
then
$$
E(K) \otimes\zp = \zp(y_K \otimes 1) \simeq \zp,\qquad \sha(E/K)[p^\infty] = 0.
$$
(2) Assume that $(K, p) = (\Q(\sqrt{-3}), 3)$; then $y_K \in 3 E(K)$.
If $3\nmid a_3$, assume, in addition, that $\rho(G_K)$ is not a cyclic group
of order four. If $y_K \not\in 3^2 E(K)$, then
$$
\Z_3 \simeq E(K) \otimes\Z_3 \supset 3 E(K) \otimes\Z_3 = \Z_3 (y_K \otimes 1),
\qquad \sha(E/K)[3^\infty] = 0.
$$

\sn
{\it Proof.\/} According to Proposition 5.26, the assumptions imply that
$E[p]$ is an absolutely irreducible $\F_p[G_K]$-module. If $u_K = 1$, the
statement follows from Proposition~6.5. It remains to consider the two fields
$K = \Q(i)$ and $K =\Q(\sqrt{-3})$, when $u_K = 2$ and $u_K = 3$, respectively.
We distinguish two separate cases.

\sn
{\bf Case 1:} $p \nmid u_K$ (equivalently, either $K = \Q(i)$ and $p > 2$, or
$K =\Q(\sqrt{-3})$ and $p > 3$).

\sn
We modify the constructions in [G] as follows. For any square-free integer
$n$ we let $H_n^\prime$ to be the compositum inside $H_n$ of the ring class
fields $H_\ell$, where $\ell$ runs through all prime numbers dividing $n$.
The Galois group $G_n := \Gal(H_n^\prime/H_1)$ is then canonically isomorphic
to $\prod_{\ell\mid n} G_\ell$, where $G_{\ell} = \Gal(H_\ell/H_1)$ is a
cyclic group of order $\#(G_\ell) = (\ell - \eta_K(\ell))/u_K$.
If, in addition, $(n,N) = 1$, we define
$y_n := \Tr_{H_n/H_n^\prime}(\varphi(x_n)) \in E(H_n^\prime)$.

One considers only square-free products $n$ of Kolyvagin primes $\ell$
satisfying [G, (3.1)-(3.2)]. For each such an $\ell$ fix a generator
$\sigma_\ell \in G_\ell$ and define $D_n := \prod_{\ell\mid n} D_\ell \in
\Z[G_n]$, where each $D_\ell$ is defined as in [G, \S{3}], except that
$\ell + 1$ is replaced by $\#(G_\ell) = (\ell + 1)/u_K$. The norm relation
[G, 3.7(1)] is replaced by $u_K \Tr_\ell(y_{\ell m}) = a_\ell \cdot y_m$
(which implies that [G, 3.6] still holds, since $p \nmid u_K$); the
congruence relation [G, 3.7(2)] does not change.

The points $P_n \in
E(H_n^\prime)$ are defined as in [G, (4.1)], except that we replace $H_n$
(denoted by $K_n$ in [G]) by $H_n^\prime$. The vanishing statement
$E(H_n)[p] = 0$ of [G, 4.3] (i.e., (C2) in Proposition~6.4) still holds,
by Proposition~6.5.

Kolyvagin's classes $c(n) \in H^1(K, E[p])$ are then
defined by $\res_{H_n^\prime/K} (c(n)) = \delta_n[P_n] \in H^1(H_n^\prime,
E[p])$ (hence $c(1) = \delta y_K$). These classes (and their images $d(n)
\in H^1(K, E)[p]$) have all the properties listed in [G, \S{6}-\S{7}]
(except that $H_n$ needs to be replaced by $H_n^\prime$).
In the formula [G, p. 246, l. 2] one needs to replace $Q_n$ by $u_K Q_n$,
but this is harmless for the argument proving the key statement [G, 6.2(2)],
since $p \nmid u_K$.

The rest of the proof goes through as in the situation considered in
Proposition~6.4.

\sn
{\bf Case 2:} $p \mid u_K$ (equivalently, $K =\Q(\sqrt{-3})$ and $p = u_K = 3$).

\sn
According to Proposition~4.14(1), there exist infinitely many primes
$q\nmid 3N$ satisfying $3 \nmid (q + 1 - a_q)$ (which is equivalent to
$3 \nmid (\eta_K(q) + 1 - a_q)$); fix once for all such a prime $q$.

Consider square-free products $n$ of primes $\ell \nmid 3Nq$ satisfying
Kolyvagin's condition [G, (3.2)] (which implies that $\eta_K(\ell) = - 1$, by
[G, (3.3)]). For each such $n$ we consider the point $y_n := \varphi(x_{qn})
\in E(H_{qn})$. The Galois group $G_n := \Gal(H_{qn}/H_q)$ is canonically
isomorphic to $\prod_{\ell\mid n} G_\ell$, and each $G_\ell$ is cyclic
of order $\ell + 1$. We define $D_n$ and $\Tr_\ell$ as in [G, \S{3}]. The
statements of [G, 3.6-3.7] and the definition of $P_n$ in [G, (4.1)] are
unchanged, except that each $H_n$ (denoted by $K_n$ in [G]) needs to be
replaced by $H_{qn}$ (so that $P_n \in E(H_{qn})$). One obtains again
classes $c(n) \in H^1(K, E[p])$ and $d(n) \in H^1(K, E)[p]$), with $c(1)
= \delta y_{K,q}$. They have all the properties listed in [G, \S{6}-\S{7}],
except that $H_n$ needs to be replaced by $H_{qn}$.

The rest of the proof goes through as in the situation considered in
Proposition~6.4, except that $y_K$ in [G, \S{9}-\S{10}] needs to be replaced
by $y_{K,q}$, and $P_\ell \in H_\ell$ in [G, \S{10}] by $P_\ell \in H_{q\ell}$.
The conclusion is that $\Sel_3(E/K) = (\Z/3\Z) \cdot \delta y_{K,q}$, which
is equivalent to $\sha(E/K)[3^\infty] = 0$ and $E(K) \otimes \Z_3 =
\Z_3(y_{K,q} \otimes 1) \simeq \Z_3$, since $E(K)[3] = 0$. In particular,
$E(K) \otimes \Z_3 = E(K)_{HP} \otimes \Z_3 \simeq \Z_3$, which implies that
$\Z_3(y_K \otimes 1) = 3 E(K) \otimes \Z_3$, by Proposition~4.14(4).

\sn
{\bf (6.8)} \ Combining Theorem~6.7 with Theorems~4.9 and 4.17, respectively,
we obtain the following results.

\proclaim{(6.9) Theorem}. Assume that $p \not= 2$, $E[p]$ is an irreducible
$\F_p[G_{\Q}]$-module and $p\nmid N \cdot a_p \cdot (a_p - 1) \cdot
(a_p - \eta_K(p))\cdot c_{\rm Tam}(E/\Q)$. If $y_K \not\in p E(K)$,
then the conclusions of Theorem~4.9 hold.

\proclaim{(6.10) Theorem}. Assume that $K = \Q(\sqrt{-3})$, $p = 3$, $E[3]$
is an irreducible $\F_3[G_{\Q}]$-module, $\rho(G_K)$ is not a cyclic group
of order four and $3\nmid a_3 \cdot (a_3 - 1) \cdot c_{\rm Tam}(E/\Q)$.
If $y_K \not\in 3^2 E(K)$, then the conclusions of Theorem~4.17 hold,
with the following modification: each $E(K_n) \otimes {\bf Z}_3$ is
generated over ${\bf Z}_3[\Gal(K_n/K)]$ by the traces to $\ki$ of the
Heegner points of conductors dividing $3^\infty q$, for any prime $q$
as in 0.23.

\bigskip
\centerline{\bf References}

\vglue12pt
\item{[B]} \ M. Bertolini, {\it Selmer groups and Heegner points in
anticyclotomic $\zp$-extensions}, Comp. Math. {\bf 99} (1995),
153--182.

\smallskip
\item{[BK]} \ S. Bloch, K. Kato, {\it $L$-functions and Tamagawa numbers
of motives}, in: The Grothendieck Festschrift I, Progress in Mathematics
{\bf 86}, Birkh\"auser, Boston, Basel, Berlin, 1990, pp. 333--400.

\smallskip
\item{[Ch]} \ B. Cha, {\it Vanishing of some cohomology groups and bounds
for the Shafarevich--Tate groups of elliptic curves}, J. of Number
Theory {\bf 111} (2005), 154--178.

\smallskip
\item{[CS1]} \ M. \c{C}iperiani, J. Stix, {\it Weil-Ch\^atelet divisible
elements in Tate-Shafarevich groups I: The Bashmakov problem for elliptic
curves over $\Q$}, Comp. Math. {\bf 149} (2013), 729--753.

\smallskip
\item{[CS2]} \ M. \c{C}iperiani, J. Stix, {\it Weil-Ch\^atelet divisible
elements in Tate-Shafarevich groups II: On a question of Cassels},
J. reine angew. Math. {\bf 700} (2015), 175--207.

\smallskip
\item{[Co1]} \ C. Cornut, {\it Reduction de Familles de points CM}, thesis,
2000.

\smallskip
\item{[Co2]} \ C. Cornut, {\it Mazur's conjecture on higher Heegner points},
Invent. Math. {\bf 148} (2002), 495--523.

\smallskip
\item{[Di]} \ L.E. Dickson, {\it Linear groups: With an exposition of the
Galois field theory}, Dover, New York, 1958.

\smallskip
\item{[FD]} \ B. Farb, R.K. Dennis, {\it Noncommutative Algebra}, 
Graduate Texts in Math. {\bf 144}, Springer, Berlin, 1993.

\smallskip
\item{[FoPR]} \ J.-M. Fontaine, B. Perrin-Riou, {\it Autour des conjectures
de Bloch et Kato: cohomologie galoisienne et valeurs de fonctions $L$}, in:
Motives (Seattle, 1991), Proc. Symp. in Pure Math. {\bf 55/I}, American
Math. Society, Providence, Rhode Island, 1994, pp. 599--706.

\smallskip
\item{[GJPST]} \ G. Grigorov, A. Jorza, S. Patrikis, W.A. Stein, C. Tarnita,
{\it Computational verification of the Birch and Swinnerton-Dyer conjecture
for individual elliptic curves}, Math. of Comp. {\bf 78} (2009),
No. 268, 2397--2425.

\smallskip
\item{[G]} \ B.H. Gross, {\it Kolyvagin's work on modular elliptic curves},
in: $L$-functions and arithmetic (Durham, 1989; J. Coates, M.J. Taylor, eds.),
LMS Lect. Note Ser. {\bf 153}, Cambridge Univ. Press, Cambridge, 1991,
pp. 235--256.

\smallskip
\item{[GZ]} \ B.H. Gross, D.B. Zagier, {\it Heegner points and derivatives
of $L$-series}, Invent. Math. {\bf 84} (1986), 225--320.

\smallskip
\item{[Gu]} \ R.M. Guralnick, {\it Small representations are completely
reducible}, J. of Algebra {\bf 220} (1999), 531--541.

\smallskip
\item{[H1]} \ B. Howard, {\it The Heegner point Kolyvagin system}, Comp.
Math. {\bf 140} (2004), 1439--1472.

\smallskip
\item{[J]} \ D. Jetchev, {\it Global divisibility of Heegner points and
Tamagawa numbers}, Comp. Math. {\bf 144} (2008), 811--826.

%\smallskip
%\item{[JSW]} \ D. Jetchev, C. Skinner, X. Wan, {\it The Birch and
%Swinnerton-Dyer formula for elliptic curves of analytic rank one},
%Cambridge J. Math. {\bf 5} (2017), 369--434.

\smallskip
\item{[K1]} \ V. A. Kolyvagin, {\it Euler systems}, in: The Grothendieck
Festschrift, vol. II, Progress in Math. {\bf 87}, Birkh\"auser, Boston,
Basel, Berlin, 1990, pp. 435--483.

\smallskip
\item{[K2]} \ V. A. Kolyvagin, {\it On the structure of Shafarevich--Tate
groups}, in: Proc. USA--USSR Symposium on Algebraic Geometry, Chicago,
1989, Lect. Notes in Math. {\bf 1479}, Springer, Berlin, 1991, pp. 94--121.

\smallskip
\item{[LW]} \ T. Lawson, C. Wuthrich, {\it Vanishing of some Galois cohomology
groups for elliptic curves}, in: Elliptic curves, modular forms and Iwasawa
theory (in honour of John H. Coates' 70th birthday, Cambridge 2015; eds.
D.~Loeffler, S.L.~Zerbes), Springer Proc. in Math. and Stat. {\bf 188},
Springer, 2016, pp. 373--399.

\smallskip
\item{[Ma]} \ A. Matar, {\it Selmer groups and Anticyclotomic
$\zp$-extensions}, Math. Proc. Camb. Phil. Soc. {\bf 161} (2016),
409--433.

\smallskip
\item{[Mz]} \ B. Mazur, {\it Modular curves and arithmetic}, in: Proc.
ICM 1983 (Warsaw), Vol. 1, PWN, Warsaw, 1984, pp. 185--211.

\smallskip
\item{[N1]} \ J. Nekov\'a\v{r}, {\it On the parity of Selmer groups II},
C. R. Acad. Sci. Paris, S\'er. I Math. {\bf 332} (2001), no. 2, 99--104.

\smallskip
\item{[N2]} \ J. Nekov\'a\v{r}, {\it Selmer complexes}, Ast\'erisque
{\bf 310} (2006), Soc. Math. de France, Paris.

\smallskip
\item{[N3]} \ J. Nekov\'a\v{r}, {\it The Euler system method for CM points
on Shimura curves}, in: $L$-functions and Galois representations (Durham,
July 2004), LMS Lect. Note Ser. {\bf 320}, Cambridge Univ. Press, 2007,
pp. 471--547.

\smallskip
\item{[PR]} \ B. Perrin-Riou, {\it Fonctions $L$ $p$-adiques, th\'eorie
d'Iwasawa et points de Heegner}, Bull. Soc. Math. France {\bf 115} (1987),
399--456.

\smallskip
\item{[Sa]} \ C.-H. Sah, {\it Automorphisms of finite groups}, J. of Algebra
{\bf 10} (1968), 47--68.

\smallskip
\item{[S1]} \ J.-P. Serre, {\it Propri\'et\'es galoisiennes des points
d'ordre fini des courbes elliptiques}, Invent. Math. {\bf 15} (1972),
259--331.

\smallskip
\item{[S2]} \ J.-P. Serre, {\it Sur la semisimplicit\'e des produits tensoriels
de repr\'esentations de groupes}, Invent. Math. {\bf 116} (1994), 513--530.

\smallskip
\item{[Va]} \ V. Vatsal, {\it Special values of anticyclotomic $L$-functions},
Duke Math. J. {\bf 116} (2003), 549--566.

\smallskip
\item{[Z]} \ S.-W. Zhang, {\it Gross--Zagier formula for $GL_2$}, Asian J.
Math. {\bf 5} (2001), 183--290.

\vskip12pt

\noindent
Ahmed Matar, Department of Mathematics, University of Bahrain,
P.O. Box 32038, Sukhair, Bahrain

\vskip12pt

\noindent
Jan Nekov\'a\v{r}, Sorbonne Universit\'e, Campus Pierre et Marie Curie,
Institut de Math\'ematiques de Jussieu, Th\'eorie des Nombres, Case 247,
4 place Jussieu, 75252 Paris cedex 05, France

\bye